\documentclass[11 pt,leqno]{amsart}

\usepackage{amsmath}
\usepackage{amssymb}
\usepackage{amsfonts}
\usepackage{amsthm}
\usepackage{hyperref}
\usepackage{enumerate}

\usepackage{bbm}
\usepackage[mathscr]{eucal}

\usepackage{tikz}
\usetikzlibrary{intersections}
\usetikzlibrary{calc}
\usepackage{bbm}
\usepackage{color}\usepackage{enumitem}
\usepackage[mathscr]{eucal}

\theoremstyle{plain}
\newtheorem{thm}{Theorem}[section]
\newtheorem{lemma}[thm]{Lemma}

\newtheorem{proposition}[thm]{Proposition}
\newtheorem{problem}[thm]{Problem}

\newtheorem{conjecture}[thm]{Conjecture}
\theoremstyle{remark}
\newtheorem{remark}[thm]{Remark}

\newtheorem{examples}[thm]{Examples}

\numberwithin{equation}{section}

\def\BlackAndWhiteOnly{0}

\newcommand{\inn}[2]{\langle#1,#2\rangle}

\newcommand{\Qfourx}[1]{\d*(\d-1)/(\d*\d+2*#1-1)}
\newcommand{\Qfoury}[1]{(\d-1)/(\d*\d+2*#1-1)}
\newcommand\Qthreex[1]{(\d-#1)/( \d-#1+1)}
\newcommand\Qthreey[1]{1/(\d-#1+1)}

\newcommand{\definecoords}{
\def\ptsize{.1pt}
\def\QbgFillOpacity{.2}
\def\QbgFillColor{black}
\def\QbgLineOpacity{.6}
\def\QbgDrawCritSeg{1}
\def\QbgCritSegStyle{solid}
\def\QbgCritSegOpacity{\QbgLineOpacity}
\def\QbbFillOpacity{.1}
\def\QbbLineOpacity{.5}
\if\BlackAndWhiteOnly1
\def\AFillColor{black}
\def\AFillOpacity{.3}
\else
\def\AFillColor{red}
\def\AFillOpacity{.1}
\fi
\def\ALineOpacity{.8}
\coordinate (Q1) at (0,0);
\coordinate (Q2) at ( {(\d-1)/(\d-1+\b)}, {(\d-1)/(\d-1+\b)} );
\coordinate (Q3) at ( {\Qthreex{\b}}, { \Qthreey{\b} } );
\coordinate (Q30) at ( {\Qthreex{0}}, { \Qthreey{0} } );
\coordinate (Q4b) at ( { \Qfourx{\b} }, { \Qfoury{\b} } );
\coordinate (Q4g) at ( { \Qfourx{\g} }, { \Qfoury{\g} } );
\coordinate (R) at ( { \d*(\d-1)/ (\d*\d-1+\b) } , { (\d-1)/(\d*\d-1+\b) } );
\coordinate (C1) at ( { (\Qfourx{\b}+\Qfourx{\g})/2 }, {(\Qfoury{\b}+\Qfoury{\g})/2} );
\coordinate (C2) at (Q4b);
}

\newcommand{\drawauxlines}[2]{
\draw (0,0) [->] -- (0,1) node [left] {$\frac1q$};
\draw (0,0) [->] -- (1,0) node [below] {$\frac1p$};
\draw[dashed,opacity=.3] (1,0) -- (0,1);
\draw[dashed,opacity=.3] (#1) -- (1,{1/\d});
\draw[dashed,opacity=.3] (#2) -- (1,1);
\draw[dashed,opacity=.1]
(.5, 0) -- (.5, .5);
}

\newcommand{\drawQbg}{
\fill (Q1) node [left] {$Q_1$} circle [radius=.02em];
\fill (Q2) node [above left] {$Q_{2,\beta}$} circle [radius=\ptsize];
\fill (Q3) node [right] {$Q_{3,\beta}$} circle [radius=\ptsize];
\fill (Q4g) node [below right] {$Q_{4,\gamma}$} circle [radius=\ptsize];
\fill[color=\QbgFillColor,opacity=\QbgFillOpacity] (Q1) -- (Q2) -- (Q3) -- (Q4g) -- cycle;
\draw[opacity=\QbgLineOpacity] (Q1) -- (Q2) -- (Q3);
\draw[opacity=\QbgLineOpacity] (Q4g) -- (Q1);
\if\QbgDrawCritSeg1
\draw[style=\QbgCritSegStyle,opacity=\QbgCritSegOpacity] (Q3) -- (Q4g);
\fi
}

\newcommand{\drawQbb}{
\fill[opacity=\QbbFillOpacity] (Q4g) -- (Q4b) -- (Q3) -- cycle;
\draw[opacity=\QbbLineOpacity] (Q4g) -- (Q4b) -- (Q3);
}

\newcommand{\lblQbb}{
\fill (Q4b) node [below right] {$Q_{4,\beta}$} circle [radius=\ptsize];
}

\newcommand{\lblQthreezero}{
\fill (Q30) node [below] {$Q_{3,0}$} circle [radius=\ptsize];
}

\newcommand{\drawA}{
\fill[color=\AFillColor,opacity=\AFillOpacity] (Q4g) .. controls (C1) and (C2) .. (Q3);
\draw[opacity=\ALineOpacity] (Q4g) .. controls (C1) and (C2) .. (Q3);
}

\newcommand{\Pthreex}[1]{(\d-1)/(\d-1+#1)}
\newcommand{\Pthreey}{0}
\newcommand{\Pfourx}[1]{(\d-1)/(\d-1+#1)}
\newcommand{\Pfoury}[1]{(1-\d)/(\d-1+#1)}

\newcommand{\Aonex}{(\d-1)/(\d)}
\newcommand{\Aoney}{(1-\d)/(\d)}
\newcommand{\Atwox}[1]{(\d-1)/(\d-1+#1)}
\newcommand{\Atwoy}[1]{(1-\d)/(\d-0.75+#1)}

\newcommand{\definecoordsmod}{
\def\ptsize{.1pt}
\def\WbgFillOpacity{.2}
\def\WbgFillColor{gray}
\def\WbgLineOpacity{.6}
\def\WbgDrawCritSeg{1}
\def\WbgCritSegStyle{solid}
\def\WbgCritSegOpacity{\QbgLineOpacity}
\def\WbbFillOpacity{.1}
\def\WbbLineOpacity{.5}
\def\AFillColor{white}
\def\AFillOpacity{0}
\def\ALineOpacity{.8}
\coordinate (P1) at (0,0);
\coordinate (P2) at (0,1);
\coordinate (P3) at ( {\Pthreex{\b}}, { \Pthreey } );
\coordinate (P4g) at ( { \Pfourx{\b} }, { \Aoney } );
\coordinate (P4) at ( { \Pfourx{\b} }, { \Pfoury{\b} } );
\coordinate (A1) at ( { \Aonex }, { \Aoney } );
\coordinate (A2) at ( { \Atwox{\b} }, { \Atwoy{\b} } );
\coordinate (A3) at (1,0);
\coordinate (A4) at (0,1);
\coordinate (C1) at ( { (\Pfourx{\b}+2*\Aonex)/3 }, {(\Pfoury{\b}+2*\Aoney)/3} );
\coordinate (C2) at ( { \Pfourx{\b} }, {(\Pfoury{\b}+\Aoney)/2} );
}

\newcommand{\drawauxlinesmod}[1]{
\draw (0,-1) [->] -- (0,1.2) node [left] {\tiny $\frac{\alpha}{(d-1)p}$};
\draw (0,0) [->] -- (1.2,0) node [below] {\tiny $\frac1p$};
\draw[dashed,opacity=.3] (0,0) -- (1,-1);
\draw[dashed,opacity=.3] (0,1) -- (#1);
\draw[dashed,opacity=.3] ({\Aonex},1) -- ({\Aonex},-1);
\draw[dashed,opacity=.3] ({\Pthreex{\b}},1) -- ({\Pthreex{\b}},-1);
\draw[dashed,opacity=.3] (0,-{\Aonex}) -- (1,-{\Aonex});
\draw[dashed,opacity=.3] (0,-{\Pthreex{\b}}) -- (1,-{\Pthreex{\b}});
}

\newcommand{\drawWbg}{
\fill ({\Aonex},0) node [below left] {\tiny $x_1$} circle [radius=0.015em];
\fill (P3) node [below right] {\tiny $x_\beta$} circle [radius=0.015em];
\fill (A3) node [below] {\tiny $1$} circle [radius=0.015em];
\fill (A4) node [left] {\tiny $1$} circle [radius=0.015em];
\fill (0,-{\Aonex}) node [left] {\tiny $-x_1$} circle [radius=0.015em];
\fill (0,-{\Pthreex{\b}}) node [left] {\tiny $-x_\beta$} circle [radius=0.015em];
\fill[color=\WbgFillColor,opacity=\WbgFillOpacity] (P1) -- (A1) .. controls (C1) and (C2) .. (A2) -- (P3) -- (P2) -- cycle;
\draw[opacity=\WbgLineOpacity] (P1) -- (A1) .. controls (C1) and (C2) .. (A2) -- (P3) -- (P2) -- cycle;
}

\newcommand{\drawAmod}{
\fill[color=\AFillColor,opacity=\AFillOpacity] (A1) .. controls (C1) and (C2) .. (A2) -- (P4) -- cycle;
\draw[opacity=\ALineOpacity] (A1) .. controls (C1) and (C2) .. (A2) -- (P4) -- cycle;
}

\begin{document}

\title{Problems on spherical maximal functions}

\author[J. Roos]{Joris Roos}
\address{Department of Mathematics and Statistics, University of Massachusetts Lowell, Lowell, MA, USA}
\email{joris\_roos@uml.edu}

\author[A. Seeger]{Andreas Seeger}
\address{Department of Mathematics, University of Wisconsin--Madison, Madison, WI, USA}
\email{aseeger@wisc.edu}

\begin{abstract} We survey old and new conjectures and results on various types of spherical maximal functions, emphasizing problems with a fractal dilation set.
\end{abstract}


\maketitle

\section{Introduction}\label{sec:intro}
For a locally integrable function $f$ on ${{\mathbb{R}}}^d$ with $d\ge 2$ we let $\sigma$ denote the rotation invariant measure on $S^{d-1}=\{x\in {{\mathbb{R}}}^d:|x|=1\}$, i.e. the normalized surface measure on $S^{d-1}$. Let
\[A_t f(x) = \int f(x-ty)\, d\sigma(y),\]
the average of $f$ over the sphere of radius $t$ centered at $x$. We may write $A_t f=f*\sigma_t$ where the dilate $\sigma_t$ is given
by $\inn{\sigma_t}{f}=\inn{\sigma}{f(t\cdot)} $.
One is interested in the pointwise behavior of
$A_t f(x)$ as $t\to 0$.
It is clear that for continuous $f$ the averages $ A_t f$ converge to $f$ as $t\to 0$, uniformly on compact sets. For larger classes (such as $L^p_{{\text{\rm loc}}}$ classes and weighted analogues) one is interested in proving almost everywhere convergence. Such results can be seen as a singular version of the classical Lebesgue differentiation result in which spheres are replaced by solid balls; many other singular variants are discussed in the 1978 survey by Stein and Wainger \cite{SteinWaingerBull}.

The proof of the a.e. convergence results requires quantitative bounds on the maximal operator $f\mapsto \sup_{t>0}|A_t f|$.
The known results are due to Stein \cite{SteinPNAS1976}, for $d\ge 3$, and to Bourgain \cite{BourgainJdA1986}, for $d=2$. They showed that the maximal operator is bounded on $L^p({{\mathbb{R}}}^d)$ if and only if
for $p>\frac{d}{d-1}$.
The case $d=2$ is harder than the higher-dimensional cases, since a crucial $L^2$ estimate fails for the maximal operator (see also \cite{MockenhauptSeegerSogge1992}).
The bounds for the maximal operator imply a.e. convergence for all functions in $L^p_{{\text{\rm loc}}} ({{\mathbb{R}}}^d)$ for $p>\frac{d}{d-1}$. This result is sharp also by probabilistic arguments in Nikishin--Stein theory (see \cite{Stein1961}, and \cite[Ch. VI]{GCRdF}, \cite[Ch. X.2]{Stein1993}).

In this survey we concentrate on the situation where the allowable radii are restricted to a given set ${{\mathcal{E}}}\subset {{\mathbb{R}}}^+=(0,\infty)$. It is then natural to ask whether for all $f\in L^p_{{\text{\rm loc}}} $ we have $A_t f(x)\to f(x)$ for almost every $x\in {{\mathbb{R}}}^d$ as $t$ tends to $0$ within ${{\mathcal{E}}}$
(for this statement we implicitly assume that $0$ is an accumulation point of ${{\mathcal{E}}}$). Such a result can be proved using quantitative bounds for the maximal function
\begin{equation}\label{eq:MGa}{{\mathcal{M}}}_{{\mathcal{E}}} f(x) =\sup_{t\in {{\mathcal{E}}}} |A_t f(x)|, \end{equation}
which is a priori well defined as a measurable function at least for continuous $f$. One needs to show that ${{\mathcal{M}}}_{{\mathcal{E}}}$ satisfies a weak type $(p,p)$ inequality, i.e.
\[{{\text{\rm meas}}}(\{x\in{{\mathbb{R}}}^d\,:\, {{\mathcal{M}}}_{{\mathcal{E}}} f(x)>{\lambda}\}){\lesssim} {\lambda}^{-p}\|f\|_p^p.\]
This means $L^p\to L^{p,\infty}$ boundedness where $L^{p,\infty}$ is the familiar Marcinkiewicz space (the largest in the Lorentz-scale $L^{p,r}$ for $0<r\le \infty$). By interpolation with the $L^\infty$ bound this implies that ${{\mathcal{M}}}_{{\mathcal{E}}}$ is bounded on $L^{\tilde p} $ for $p<\tilde p\le \infty$.
We are also interested in weighted analogues, where for suitable positive weight functions $w$ the weak type estimate takes the form
\begin{equation} \label{eq:weaktypeMGa} \int_{\{{{\mathcal{M}}}_{{\mathcal{E}}} f(x)>{\lambda}\}} w(x)\,dx {\lesssim} {\lambda}^{-p} \|f\|_{L^p(w)}^p, \end{equation}
and $\|f\|_{L^p(w)}=(\int|f(x)|^p w(x)\,dx)^{1/p}$ is a weighted $L^p$ norm.

\subsubsection*{\texorpdfstring{\it Structure of this survey.}{Structure of this survey}} In \S\ref{sec:lacunary} we consider lacunary dilation sets and associated problems on the lacunary spherical maximal function, in particular their behavior in spaces near $L^1$. In \S\ref{sec:ptwise} we cover the $L^p\to L^p$ and $L^p\to L^{p,\infty}$ boundedness problem for the operators ${{\mathcal{M}}}_{{\mathcal{E}}}$, and formulate conjectures for various endpoint cases. Weighted analogues and the relation to sparse domination are briefly discussed in \S\ref{sec:sparse}. The sparse domination problem for ${{\mathcal{M}}}_{{\mathcal{E}}}$ is closely related to {\em $L^p$ improving bounds} for a local maximal operator ${{\mathcal{M}}}_E$ with dilation set $E\subset [1,2]$ in \S\ref{sec:Lpimpr}. There we discuss known results involving sufficient conditions and necessary conditions and state
the main conjecture on $L^p\to L^q$ inequalities for ${{\mathcal{M}}}_E$. This conjecture involves the Legendre--Assouad function $\nu_{{\mathcal{E}}}^\sharp$ which plays a central role in several problems on operators in Fourier analysis associated with restricted dilation sets. Another such instance, the problem of weighted norm inequalities for ${{\mathcal{M}}}_{{\mathcal{E}}}$ in spaces with power weights, is covered in \S\ref{sec:power}.
In \S\ref{sec:further} we discuss further directions on multiparameter analogues, on helical maximal functions and Nevo--Thangavelu maximal functions on step two nilpotent Lie groups. We conclude in \S\ref{sec:lsm} with a fractal $L^p$ local smoothing conjecture for the wave equation. We also state a related $L^p\to L^q$ version of this conjecture which is closely related to the above mentioned problem on $L^p$ improving properties for the local maximal operators. Again the Legendre--Assouad function plays a central role.

\section{Lacunary dilation sets}\label{sec:lacunary} The lacunary case, e.g. ${{\mathcal{E}}}^{\mathrm{lac}}=\{2^k:k\in {{\mathbb{Z}}}\}$
was considered by
C. \!\!P. Calder\'on \cite{CalderonCP}
who proved that ${{\mathcal{M}}}_{{{\mathcal{E}}}^{\mathrm{lac}}}$ is bounded on $L^p$ for $p>1$ (see also Coifman--Weiss \cite{CoifmanWeissBookReview}).
In order to use Calder\'on--Zygmund theory one employs the classical dyadic frequency decomposition
\[\sigma= \sum_{j=0}^\infty \sigma_j,\]
where the Fourier transform $\widehat \sigma_j(\xi) $ is supported on the annulus $\{|\xi|\approx 2^j\}$. One then considers for each $j$ the maximal function
\[M^j f= \sup_{k\in\mathbb{Z}} 2^{-kd}|\sigma_j(2^{-k} \cdot)*f|.\]
It satisfies an $L^2$-bound \[\|M^j\|_{L^2\to L^2} = O(2^{-j(d-1)/2})\] which follows from the stationary phase bound
\begin{equation*}
\|\widehat \sigma_j\|_\infty= O(2^{-j(d-1)/2} )\end{equation*} and Plancherel's theorem. Indeed,
\begin{align*} \|M^jf\|_2 &\le \Big\| \Big(\sum_{k\in{{\mathbb{Z}}}}|2^{-kd}\sigma_j(2^{-k}\cdot)*f|^2\Big)^{\frac 12} \Big\|_2 = \Big(\sum_{k\in{{\mathbb{Z}}}} \big\| \widehat \sigma_j(2^k\cdot) \widehat f\big\|_2^2\Big)^{1/2}\\& {\lesssim} 2^{-j\frac{d-1}{2}} \|\widehat f\|_2=2^{-j\frac{d-1}{2}} \| f\|_2.
\end{align*}
For $p=1$ the operator $M^j$
satisfies a weak type $(1,1)$ inequality with constant $O(1+j)$ as $j\to \infty$. This can be shown by arguments from a Calder\'on--Zygmund theory.
Then the Marcinkiewicz interpolation theorem is used to get
$\|M^j\|_{L^p\to L^p} =O(2^{-j{\varepsilon}(p)}) $ with ${\varepsilon}(p)>0$ for $1<p<\infty$, with the implicit constant depending on $p$.
This interpolation argument only provides limited information at the endpoint $p=1$ (except for functions locally in $L\log L$, see \cite{ChristStein}).

\begin{problem} \label{pr:lac} For ${{\mathcal{E}}}^{\mathrm{lac} }= \{2^k:k\in {{\mathbb{Z}}}\}$, is ${{\mathcal{M}}}_{{{\mathcal{E}}}^{\mathrm{lac}}}$ of weak type $(1,1)$?
\end{problem}

Christ \cite{Christ-rough} showed by more sophisticated arguments that a weak type $(1,1)$ estimate for ${{\mathcal{M}}}_{{{\mathcal{E}}}^{\mathrm{lac}}}$ holds for functions in the Hardy space $H^1$.
The proof relies on the atomic decomposition of $H^1$
(see e.g. \cite{CoifmanWeissatomic})
which breaks down for the space $L^1$. We note that a more standard $L^1\to L^{1}$ estimate fails for the maximal operator. Since the Marcinkiewicz space is not normable one cannot reduce to estimates on single atoms as in the standard theory of singular integral operators. Instead one proves $L^2$ estimates in the complement of a suitable exceptional set which is constructed based on the size of the coefficients in the atomic decomposition.

The $H^1\to L^{1,\infty}$ result of Christ can be extended to a much larger class of lacunary maximal functions generated by compactly supported measures, see \cite{SW-Marc}. One imposes a dimensional assumption (that the support of the measure can be covered with $O(\delta^{-a})$ balls of radius $\delta$) together with a corresponding sharp $L^p\to L^{p}_{a/p'}$ Sobolev type estimate for some $p\in (1,2]$, for functions appropriately localized in frequency. Thus instead of having to prove $L^2$ estimates off an exceptional set it suffices to show certain $L^p$ bounds off that set. This method applies to lacunary maximal functions for averages over hypersurfaces where the Gaussian curvature does not vanish to infinite order \cite{SoggeStein},
and also for averages over finite type curves in ${{\mathbb{R}}}^3$ \cite{PramanikSeeger} and in higher dimensions \cite{KoLeeOh2023}. We emphasize that here we are considering lacunary maximal operators generated by standard {\it isotropic} dilations.

\begin{remark}\label{rem:Christ-2par}
Christ \cite{Christ-rough} also showed a more sophisticated result regarding a maximal function for averages {\it along the curve} $\gamma(t)=(t,t^2)$, namely
\[\sup_{k\in{{\mathbb{Z}}}} \Big|2^k \int_0^{2^{-k} } f(x_1-s,x_2-s^2) ds\Big|.\] Here the underlying dilation structure is nonisotropic, with the parabolic dilations given by $x\mapsto (t x_1, t^2 x_2)$, where $t=2^{-k}$. Christ
proved an $H^1_{\mathrm{par}} ({{\mathbb{R}}}^2)\to L^{1,\infty} ({{\mathbb{R}}}^2)$ bound for the maximal operator, where $H^1_{\mathrm{par}}$ is the Hardy space associated with the parabolic dilations. It is conjectured that similar results hold for the moment curves in higher dimensions, with the appropriate nonisotropic Hardy-space, but this is still open.
As in Problem \ref{pr:lac} one can also conjecture the more elusive weak type $(1,1)$ inequality for all $f\in L^1$.
\end{remark}

The ideas by Christ have been refined to make progress towards the weak type $(1,1)$ inequality in Problem \ref{pr:lac}. Indeed, Tao, Wright and the second author \cite{STW-loglog} obtained
a weak type $L\log\log L$ inequality for the lacunary spherical maximal operator, i.e.
\[{{\text{\rm meas}}}\big( \{ x: {{\mathcal{M}}}_{{{\mathcal{E}}}^{\mathrm{lac} }}f(x)>\alpha\}\big) {\lesssim} \int \frac{ |f(x)|}{\alpha} \log\log \big(10+ \frac{|f(x)|}{\alpha} \big) \, dx . \]
This implies almost everywhere convergence results for functions locally in $L\log\log L$.
In the proof one still works with $L^2$ bounds off an exceptional set. However instead of using an atomic decomposition one constructs an exceptional set based on a geometric decomposition using notions of `length' and `thickness'.
Building on this decomposition, Cladek and Krause \cite{CladekKrause} provided further refinements, and
obtained a.e. convergence for functions locally in
$L \log\log\log L(\log\log\log\log L)^{1+{\varepsilon}}$, for all ${\varepsilon}>0$.

\section{\texorpdfstring{General dilation sets: ${L}^{{p}}$ bounds}{General dilation sets}} \label{sec:ptwise} We now consider
general dilation sets ${{\mathcal{E}}}\subset {{\mathbb{R}}}^+$. Our goal is to understand how the $L^p$ bounds for the spherical maximal operator ${{\mathcal{M}}}_{{\mathcal{E}}}$ depend on the set ${{\mathcal{E}}}$.
One should expect results that are intermediate between those for the lacunary case $\mathcal{E}^{\mathrm{lac}}$
and the unrestricted case $\mathcal{E}={{\mathbb{R}}}^+$.
It is natural to test intermediate dilation sets such as $\{n^{-\alpha}:n\in {{\mathbb{N}}}\}$ for $\alpha>0$, or $\{2^{-\sqrt{n}}: n\in {{\mathbb{N}}} \}$.
However, it turns out $L^p$ bounds for these sets hold in the same range as for the unrestricted case. This is due to the scaling invariance of the operator: $\|{{\mathcal{M}}}_{R{{\mathcal{E}}}}\|_{L^p\to L^p} =
\|{{\mathcal{M}}}_{{\mathcal{E}}}\|_{L^p\to L^p} $.
This makes it natural to endow ${{\mathbb{R}}}^+$ with the metric
\begin{equation}\label{eqn:metric}
d_\times(s,t) = |\log_2 (s/t)|.
\end{equation}
For $J\subset {{\mathbb{R}}}^+$, write $|J|_{\times}$ for its diameter with respect to $d_\times$. Observe that each dilation $t\mapsto \lambda t$, $\lambda>0$, acts isometrically on ${{\mathbb{R}}}^+$; in particular, $|[R,2R]|_{\times}=1$ for every $R>0$.
If $E\subset {{\mathbb{R}}}^+$ has finite diameter, we define $N(E,\delta)$ to be its {\it entropy number}: for $\delta>0$, this is the least number of intervals of diameter $\delta$ needed to cover $E$, with diameter understood throughout in the metric $d_\times$.

It was proved in \cite{SeegerWaingerWright1995} that the $L^p$ boundedness range of ${{\mathcal{M}}}_{{\mathcal{E}}}$ depends on the following scale-invariant extension of the (upper) Minkowski dimension.
For fixed ${{\mathcal{E}}}\subset {{\mathbb{R}}}^+$ define
\begin{equation}\label{eq:betaE} \beta = \beta_{{\mathcal{E}}} = \limsup_{\delta\to 0} \frac{\sup_{|J|_\times= 1} \log N({{\mathcal{E}}}\cap J, \delta)}{\log (\delta^{-1}) },
\end{equation}
where the supremum ranges over all intervals $J\subset {{\mathbb{R}}}^+$ with $d_\times$-diameter equal to one.
One can equivalently express this as a condition on the rescaled sets
\begin{equation}\label{eq:GaR}
{{\mathcal{E}}}_R= \tfrac 1R \big({{\mathcal{E}}}\cap[R, 2R]\big) \subset [1,2]
\end{equation}
with $R>0$. Then $\beta$ is the infimum over all $b>0$ such that there exists a constant $C\in (0,\infty)$ such that
\[N({{\mathcal{E}}}_R,\delta)\le C \delta^{-b}\] for all $\delta \in (0,1)$ and all $R>0$.
Note that if ${{\mathcal{E}}}$ has bounded diameter, then $\beta_{{\mathcal{E}}}$ is equal to the upper Minkowski dimension.

\begin{thm}[\cite{SeegerWaingerWright1995}] \label{thm:sww}
Let ${{\mathcal{E}}}\subset {{\mathbb{R}}}^+$, $\beta=\beta_{{\mathcal{E}}}$ and $p_{\mathrm{cr}}(\beta) = 1+\frac{\beta}{d-1}$. Then ${{\mathcal{M}}}_{{\mathcal{E}}}$ is bounded on $L^p$ for $p>p_{\mathrm{cr}}(\beta) $ and unbounded on $L^p$ for $p<p_{\mathrm{cr}}(\beta)$.
\end{thm}

\subsubsection*{\texorpdfstring{\it Examples}{Examples}}
For $E\subset [1,2]$ let
\begin{equation}\label{eq:GammaE}
{{\mathcal{E}}}
=\bigcup_{k\in {{\mathbb{Z}}}} 2^k E.
\end{equation}
Then $\beta_{{\mathcal{E}}}$ equals the classical upper Minkowski dimension of $E$. With ${{\mathcal{E}}}_R$ as in \eqref{eq:GaR}
we have
$\beta_{{{\mathcal{E}}}} \ge \sup_{R>0} \beta_{{{\mathcal{E}}}_R}$
and the inequality can be strict.

Interesting examples include the convex sequences $E = \{1+n^{-\alpha}\,:\,n\in{{\mathbb{N}}} \}$ for $\alpha>0$, where $\beta=\frac1{1+\alpha}$ and $E$ being the middle-third Cantor set in $[1,2]$, where $\beta=\log_32$.
In both cases the corresponding maximal operators ${{\mathcal{M}}}_E$
and ${{\mathcal{M}}}_{{\mathcal{E}}}$ are bounded in an intermediate range with the critical exponent given by $p_{\mathrm{cr}}(\beta)=1+\frac{\beta}{d-1}$.
Note that the convex sequences have Hausdorff dimension zero, illustrating that Hausdorff dimension does not play a role in this problem.

\subsection*{\texorpdfstring{\it A necessary condition}{A necessary condition}}
Unboundedness for $p\le \frac{d}{d-1}$ for the full maximal function was observed by Stein \cite{SteinPNAS1976}: one tests $A_{t(x) }f(x)$ for $t(x)=|x|$ and
\begin{equation}
\label{eq:stein-example} f(x)= |x|^{1-d} (\log |x|)^{-1} {{\mathbbm 1}}_{|x|<1/2}.
\end{equation}
To see the sharpness of Theorem \ref{thm:sww} one can test the rescaled operator ${{\mathcal{M}}}_{{{\mathcal{E}}}_R}$ for every $R>0$ on characteristic functions of balls of small radius $\delta>0$. The value of the maximal function is then ${\gtrsim} \delta^{d-1}$ in a union of $\approx N({{\mathcal{E}}}_R, \delta)$ annuli of width $\approx \delta$ and radius $\approx 1$. This yields the lower bound
\begin{equation} \label{eq:MGa-lowerbound}\|{{\mathcal{M}}}_{{\mathcal{E}}}\|_{L^p\to L^p} {\gtrsim} \sup_{R>0} \sup_{\delta\in (0,1)}
\delta^{ (d-1)(1-\frac 1p)} N({{\mathcal{E}}}_R, \delta)^{\frac 1p}\end{equation}
and gives unboundedness for $p<p_\mathrm{cr}(\beta)$.

\subsection*{\texorpdfstring{\it Sketch of proof of upper bounds}{Sketch of proof of upper bounds}}
In order to prove the $L^p$-boundedness for $p>p_{\mathrm{cr}}(\beta)$ one uses the dyadic frequency decomposition $\sigma=\sum_{j\ge 0} \sigma_j$ as in \S\ref{sec:lacunary}
and derives for $E={{\mathcal{E}}}_R$ an estimate for the maximal function
\begin{equation}\label{eq:Mj-def} M_E^j f(x)= \sup_{t\in E} |\sigma_{j,t} * f|\end{equation} where $\sigma_{j,t} =t^{-d}\sigma_j(t^{-1}\cdot) $.
In view of the support of $\widehat{\sigma_j}$ one can use an averaging argument to reduce to a uniform estimate for ${{\mathcal{M}}}_{E(j)}$ where $E(j)$ is an arbitrary maximal $2^{-j}$-separated subset of $E$. Note $\#E(j)\approx N(E, 2^{-j})$. One then proves the estimate
\begin{equation}\label{eq:Mj-est} \big \|M^j_{E(j)} f\big \|_p{\lesssim} 2^{-j (d-1)(1-\frac 1p)} N(E, 2^{-j})^{\frac 1p}\|f\|_p.\end{equation}
This bound follows by interpolation from the cases $p=1$ and $p=2$. For the case $p=1$ one estimates the supremum over $t\in E(j)$ by a sum. For $p=2$ one estimates this supremum by an $\ell^2$-norm over $E(j)$, then uses Plancherel's theorem together with the bound $\|\widehat{\sigma_j}\|_\infty =O(2^{-j\frac{d-1}{2}})$.

Apply the bound $N(E, \delta){\lesssim} C_{\varepsilon} \delta^{-\beta-{\varepsilon}}$ in \eqref{eq:Mj-est} and sum in $j$. Assuming $d\ge 3$ or
$\beta<1$, this yields the $L^p$ boundedness of ${{\mathcal{M}}}_E$ for $p>1+\frac{\beta}{d-1}$.
For the global maximal operator ${{\mathcal{M}}}_{{\mathcal{E}}}$ with general ${{\mathcal{E}}}\subset {{\mathbb{R}}}^+$
one combines this with arguments from Calder\'on-Zygmund theory \cite{SeegerWaingerWright1995}.

\subsection*{\texorpdfstring{{\it What happens for $p=p_{\mathrm{cr}}(\beta)$?}}{What happens at the critical exponent?}}
Since the definition of $\beta_{{\mathcal{E}}}$ does not take into account
factors that are $O(\delta^{{\varepsilon}})$ for all ${\varepsilon}>0$ it is not the right quantity to determine outcomes at this endpoint. We compare
\eqref{eq:MGa-lowerbound} and \eqref{eq:Mj-est}.
Simply taking the $j$-sum in \eqref{eq:Mj-est} is not efficient anymore. In fact, it is shown in
\cite{SeegerWaingerWright1995} that for $E\subset [1,2]$ and $p\in (1,2)$,
\begin{equation} \label{eq:p'bd}\|{{\mathcal{M}}}_E\|_{L^p\to L^p} {\lesssim} \Big(\sum_{j\ge 0}\big[ 2^{-j (d-1)(1-\frac 1p)} N(E, 2^{-j})^{\frac 1p} ]^{p'} \Big)^{1/p'} .\end{equation}
However, the known lower bounds and a
study of the maximal operators ${{\mathcal{M}}}_E$ and ${{\mathcal{M}}}_{{\mathcal{E}}}$ when acting on radial functions (\cite{SeegerWaingerWright1997}) suggest a further improvement, namely that in \eqref{eq:p'bd} the $\ell^{p'}$ norm in $j$ may be replaced by a supremum in $j$.

\begin{problem} \label{pr:MEsharp}
Let $p\in (1,\frac{d}{d-1}).$
\begin{itemize}
\item[(i)] (Local problem)
For $E\subset [1,2]$, is it true that \[ \|{{\mathcal{M}}}_E\|_{L^p\to L^p} \approx \sup_{\delta\in (0,1)}
\delta^{ (d-1)(1-\frac 1p)} N(E, \delta)^{\frac 1p}\,\,\text{?} \]
\item[(ii)] (Global problem)
For
${{\mathcal{E}}}\subset {{\mathbb{R}}}^+$, is it true that
\[ \|{{\mathcal{M}}}_{{\mathcal{E}}}\|_{L^p\to L^{p,\infty}} \approx \sup_{R>0} \sup_{\delta\in (0,1)}
\delta^{ (d-1)(1-\frac 1p)} N({{\mathcal{E}}}_R, \delta)^{\frac 1p}\,\,\text{?} \]
\end{itemize}
\end{problem}
Note that in the local problem we have a strong type $(p,p)$ bound while in the global problem we need a weak type $(p,p)$ bound.
We conjecture an affirmative answer in both problems.
As noted in \cite{SeegerWaingerWright1995}, the lower bounds in Problem \ref{pr:MEsharp} are true and based on radial counterexamples ({\emph{cf.\ }} \eqref{eq:MGa-lowerbound}). Moreover, the answer to both questions is yes if we restrict to radial $L^p$ functions \cite{SeegerWaingerWright1997}.

To gather more evidence for the conjecture to hold we consider examples of sets in $[1,2]$ with
\begin{equation}\label{eq:endpoint-cond} N(E,\delta) \approx \delta^{-\beta}.
\end{equation}
For $E=\{1+n^{-a}: n\in {{\mathbb{N}}}\}$, \eqref{eq:endpoint-cond} holds with $\beta=\frac1{a+1}$. For this sequence example and the critical exponent $p=p_{\mathrm{cr}}(\beta)$
the endpoint $L^p\to L^p$ bound for ${{\mathcal{M}}}_E$ and the endpoint weak type $(p,p)$ bound for ${{\mathcal{M}}}_{{{\mathcal{E}}}} $ (with ${{\mathcal{E}}}=\cup_k2^kE$) have been established in \cite{STW-jussieu}.
This was done using a rather sophisticated analysis based on atomic decompositions for $L^p$ functions when $p>1$.

In contrast, the endpoint inequalities for $p=p_\mathrm{cr}(\beta)$ are still open for most examples satisfying \eqref{eq:endpoint-cond}. In particular we mention a model problem for self-similar sets of dilations:
\begin{problem} Let $E$ be the Cantor middle third set in $[1,2]$ and consider the global extension ${{\mathcal{E}}}=\cup_{k\in {{\mathbb{Z}}}} 2^kE$. Let
$p_{\mathrm{cr}}= 1+\frac{\log_3 2}{d-1}$.

(i) Is ${{\mathcal{M}}}_E$ of weak type $(p_\mathrm{cr},p_\mathrm{cr})$, or even of strong type $(p_\mathrm{cr},p_\mathrm{cr})$?

(ii) Is ${{\mathcal{M}}}_{{\mathcal{E}}}$ of weak type $(p_\mathrm{cr},p_\mathrm{cr})$?
\end{problem}

In Problem \ref{pr:MEsharp} we needed to exclude the endpoint case $p=\frac{d}{d-1}$. In this case the known and the conjectured outcomes feature additional logarithmic factors:
\begin{problem} \label{pr:MEsharppd}
Let $p=\frac{d}{d-1} $.
\begin{itemize}
\item[(i)] (Local problem) For $E\subset [1,2]$, is it true that
\[ \|{{\mathcal{M}}}_E\|_{L^{p}\to L^{p} } \approx \sup_{\delta\in (0,\frac12)} \delta^{(d-1)(1-\frac1{p} )} (\log\tfrac 1\delta)^{\frac{1}{d}}N(E,\delta)^{\frac1p}\,\,\text{?} \]
\item[(ii)] (Global problem) For ${{\mathcal{E}}}\subset {{\mathbb{R}}}^+$, is it true that
\[ \|{{\mathcal{M}}}_{{\mathcal{E}}}\|_{L^{p} \to L^{p,\infty}} \approx \sup_{R>0} \sup_{\delta\in (0,\frac12)} \delta^{(d-1)(1-\frac1p)} (\log\tfrac 1\delta)^{\frac{1}{d}}N({{\mathcal{E}}}_R,\delta)^{\frac1p}\,\,\text{?}\]
\end{itemize}
\end{problem}
Again, both the lower bounds and the upper bounds on radial functions have been verified in \cite{SeegerWaingerWright1997}, which leads to the conjecture of an affirmative answer to the questions in Problem \ref{pr:MEsharppd}.

\begin{remark} Other forms of the global problems,
on characterizations of $L^p\to L^{p,r}$ boundedness of ${{\mathcal{M}}}_{{\mathcal{E}}}$ for $p\le r<\infty$, can be found in \cite{SeegerWaingerWright1997}.
However these statements with $r<\infty$ are not relevant for pointwise convergence.
\end{remark}

\begin{remark} For the full spherical maximal operator a much easier restricted weak type (i.e. $L^{p,1}\to L^{p,\infty} $) estimate for $p=\frac{d}{d-1} $ was established by Bourgain \cite{Bourgain-CompteRendu1985} in dimension $d\ge 3$. This amounts to a weak type $(p,p)$ estimate for characteristic functions of sets with finite measure. In dimension $d=2$ this restricted weak type $(2,2)$ inequality holds on radial functions \cite{Leckband}, but fails for general functions in $L^{2,1}({{\mathbb{R}}}^2)$ due to a Besicovitch set example \cite{STW-jussieu}.
For sets $E$ satisfying \eqref{eq:endpoint-cond} with $\beta<1$, a variant of Bourgain's argument gives the restricted weak type $(p,p)$ estimate for $p=p_{\mathrm{cr}}(\beta)$.
\end{remark}

\begin{remark}The proofs of
the known
results (and most of the results in this survey) strongly rely on Fourier analysis.
It is also interesting
to find purely geometric-combinatorial arguments.
This has been accomplished in some cases, see e.g. the works by Schlag \cite{Schlag1997,Schlag-Duke1998} and recently Hickman--Jan\v car \cite{HickmanJancar}.
\end{remark}

\section{Weighted norm inequalities and sparse domination} \label{sec:sparse}
In this section we are interested in weighted $L^p$ estimates for ${{\mathcal{M}}}_{{\mathcal{E}}}$, which are closely tied to sparse domination.

\begin{problem} For given ${{\mathcal{E}}}\subset {{\mathbb{R}}}^+$, find characterizations of weights $w$ for which ${{\mathcal{M}}}_{{\mathcal{E}}}$ is bounded on $L^p(w)$.
\end{problem}

Almost conclusive results are known
for inequalities with power weights
$w_\alpha(x)=|x|^{\alpha}$
(see \S\ref{sec:power} below). However, a complete characterization of weights for which $L^p(w)$-boundedness of ${{\mathcal{M}}}_{{\mathcal{E}}}$ holds
seems currently out of reach, for any ${{\mathcal{E}}}$.
We will now briefly discuss the powerful tool of {\it sparse domination} which generates weighted norm inequalities but is also interesting in its own right.

Sparse domination, in its normed and pointwise versions, was first introduced because of its
applicability to weighted norm inequalities for Calder\'on--Zygmund operators \cite{LeCZ, LeA2, CAR2014,lerner-nazarov,Lac2015,LeNew}. Here we are using the version of
{\em bilinear sparse domination} from the works of Bernicot, Frey and Petermichl \cite{bernicot-frey-petermichl} and Culiuc, Di Plinio and Ou \cite{culiuc-diplinio-ou} which extended the scope of the theory to classes of operators beyond the Calder\'on-Zygmund theory. The concept was applied to spherical maximal operators by Lacey \cite{laceyJdA19}. For a general theory with many applications and more history see
\cite{BeltranRoosSeeger-sparse} and \cite{Conde-AlonsoDiPlinioParissisVempati} and the references in those works.

We review basic definitions.
Fix a dyadic lattice ${{\mathfrak{Q}}}$ of dyadic cubes in the sense of Lerner and Nazarov \cite{lerner-nazarov} so that the cubes of sidelength $2^{-k}$ form a grid and every cube in ${{\mathfrak{Q}}}$ of sidelength $2^{-k}$ is contained in a unique parent cube $Q'\in {{\mathfrak{Q}}}$ of sidelength $2^{1-k}$. Moreover, every compact set is contained in some $Q\in {{\mathfrak{Q}}}$.
A collection ${{\mathfrak{S}}}\subset {{\mathfrak{Q}}}$ is called {\it sparse} if for every $Q\in {{\mathfrak{S}}}$ there is a measurable subset $E_Q\subset Q$ so that $|E_Q|>\tfrac 12 |Q|$ and the sets $E_Q$ with $Q\in {{\mathfrak{S}}}$ are pairwise disjoint (the choice of constant $\tfrac 12$ is not essential here). For $f\in L^1_{{\text{\rm loc}}}$ and $Q\in {{\mathfrak{Q}}}$ we set ${{\mathrm{av}}}_Qf=\frac{1}{|Q|}\int_Qf(x)dx$. Given $1\le p_1,p_2\le \infty$ and a sparse family ${{\mathfrak{S}}}$ of dyadic cubes, the corresponding {\it sparse form} is defined by
\[{\Lambda}_{p_1,p_2}^{{\mathfrak{S}}} (f,g) = \sum_{Q\in {{\mathfrak{S}}}}|Q|\, \big({{\mathrm{av}}}_Q|f|^{p_1} \big)^{1/p_1} \big({{\mathrm{av}}}_Q|g|^{p_2} \big)^{1/p_2}.
\]

The maximal $(p_1,p_2)$-form ${\Lambda}^*_{p_1,p_2} $ is given by
\[{\Lambda}^*_{p_1,p_2} (f,g)=\sup_{{{\mathfrak{S}}}\text{ sparse}}{\Lambda}^{{\mathfrak{S}}}_{p_1,p_2} (f,g)\]
where the supremum is taken over all sparse families ${{\mathfrak{S}}}\subset{{\mathfrak{Q}}}$. We say that a sublinear operator $T$
satisfies a $(p_1,p_2)$-sparse domination inequality if for all simple functions $f$, $g$ the inequality
\begin{equation}\label{eq:sparse-dom} \Big|\int Tf(x) g(x) \,dx\Big| \le C {\Lambda}^*_{p_1,p_2}(f,g) \end{equation} holds with some constant $C$ independent of $f$ and $g$.
In this case we write $T\in \text{Sp}(p_1,p_2)$ and the best constant in this inequality is denoted by $\|T\|_{{{\mathrm{Sp}}}(p_1,p_2)}$.

Note that \begin{equation}
\label{eq:estmaxform} {\Lambda}^*_{p_1,p_2} (f,g) \le \int \big(M[|f|^{p_1} ](x) \big)^{1/p_1} \big(M[|g|^{p_2} ](x)\big) ^{1/p_2}\,dx,
\end{equation}
so the case $p_2< p_1'$ is of interest for $L^{p_1}$-bounded operators $T$.

Bernicot, Frey and Petermichl \cite{bernicot-frey-petermichl} showed that sparse domination inequalities of the form \eqref{eq:sparse-dom} yield certain weighted norm inequalities involving Muckenhoupt and reverse H\"older conditions on the weight. Recall that the Muckenhoupt $A_\rho$ characteristic of a nonnegative weight $w$ is defined by \[[w]_{A_\rho } =\sup_Q {{\mathrm{av}}}_Q[w] \big({{\mathrm{av}}}_Q [w^{-\frac{1}{\rho-1}}] \big)^{\rho-1}\] (at least for $1<\rho<\infty$) and $A_\rho$ is the class of weights with finite $A_\rho$-characteristics. For $\sigma>1$ the reverse H\"older class $\mathrm{RH}_\sigma$ consists of weights satisfying
$\big ({{\mathrm{av}}}_Q [w^\sigma]\big)^{1/\sigma }\le C {{\mathrm{av}}}_Q w$
for all cubes $Q$, and the best constant $C$ is called the reverse H\"older characteristic $[w]_{\mathrm{RH}_\sigma}$. Assuming that $p_1<p<p_2$, $p_2\le p_1'$, and $T \in {{\mathrm{Sp}}}(p_1,p_2)$ it is shown in \cite{bernicot-frey-petermichl} that under the assumption
$w\in A_{ {p}/{p_1} } \cap \mathrm{RH}_{ ({p_2'}/{p})' }$
the weighted norm inequality
\begin{equation} \label{eq:weightedineq}\|Tf\|_{L^p(w)} \le C(p_1,p_2,p, w) \|T\|_{{{\mathrm{Sp}}}(p_1,p_2)} \|f\|_{L^p(w)}
\end{equation}
holds. A major motivation for the theory of sparse domination inequalities was that it yields efficient descriptions on how the constants in the weighted norm inequality depend on the Muckenhoupt and reverse H\"older characteristics; indeed \cite{bernicot-frey-petermichl} specifies that
\[C(p_1,p_2,p, w){\lesssim} ( [w]_{A_{p/p_1}} [w]_{\mathrm{RH}_{ (p_2'/p)'}} )^\upsilon\]
with $\upsilon= \max\big\{\frac 1{p-p_1}, \frac{p_2'-1}{p_2'-p} \big\}.$

Lacey \cite{laceyJdA19} proved close to sharp sparse domination inequalities for the full and also the lacunary spherical maximal operators (and via the above reduction deduced suitable weighted norm inequalities). He established a connection with certain $L^p\to L^q$ results for a local spherical maximal function $\sup_{1\le t\le 2} |A_t f(x)|$. Essentially sharp $L^p$ improving estimates are due to Schlag \cite{Schlag1997} and Schlag and Sogge \cite{SchlagSogge1997}.
For the full spherical maximal operator these give
$L^p\to L^q$ boundedness for $(\frac 1p,\frac 1q)$ in the interior of the quadrangle with corners $Q_1=(0,0)$, $Q_2=(\frac{d-1}{d}, \frac{d-1}{d})$, $Q_3=(\frac{d-1}{d}, \frac 1{d-1})$ and $Q_4=(\frac{d(d-1)}{d^2+1},\frac{d-1}{d^2+1})$.
Most endpoint inequalities on the maximal function were settled by Lee \cite{Lee03}; see also the more recent work \cite{BeltranOberlinRoncalSeegerStovall} on $L^p$-improving estimates for variation norm operators associated to the spherical means.

Using the $L^p\to L^q$ bounds Lacey
\cite{laceyJdA19} obtained a close to sharp sparse domination result for the full spherical maximal operator ${{\mathcal{M}}}$, namely, ${{\mathcal{M}}}\in {{\mathrm{Sp}}}(p_1,p_2)$ if $(p_1^{-1}, p_2^{-1})$ is in the interior of the quadrangle with corners
\begin{equation}\label{eq:laceysp} P_1= (0,1), \,\, P_2= \big(\tfrac{d-1}{d},\tfrac{1} d\big), \,\, P_3= \big (\tfrac{d-1}{d}, \tfrac{d-1}{d} \big), \,\, P_4=\big (\tfrac{d^2-d}{d^2+1}, \tfrac{d^2-d+2}{d^2+1} \big).
\end{equation}
Using the
$L^p$-improving results for convolutions of spherical measure \cite{StrichartzTAMS1970, Littman},
he also obtained the interior of the sparse exponent set for the lacunary maximal operator, namely the interior of the triangle with corners $(0,1)$, $(1,0)$ and $(\frac{d}{d+1},\frac d{d+1})$.
Moreover, Lacey suggested the study of sparse domination for variation norm operators associated with the spherical means (which led to the paper \cite{BeltranOberlinRoncalSeegerStovall}) and also the study of sparse domination for the spherical maximal functions with restricted dilation sets.

The connection of sparse domination results with $L^p$ improving results for local rescaled operators has been noted in the literature for many other operators. A very general theorem was proved in our joint memoir with David Beltran \cite{BeltranRoosSeeger-sparse} which contains more historical references.
We introduce terminology to state this connection for the special case of the operators ${{\mathcal{M}}}_{{\mathcal{E}}}$. With ${{\mathcal{E}}}_R$ as in \eqref{eq:GaR} we let
the {\it Lebesgue exponent set}
$\mathscr{L}({{\mathcal{E}}})$ consist of all $(\tfrac1p,\tfrac1q)$ for which the three quantities
$\|{{\mathcal{M}}}_{{\mathcal{E}}}\|_{L^p\to L^{p,\infty} }$, $ \|{{\mathcal{M}}}_{{\mathcal{E}}}\|_{L^{q,1}\to L^q}$, and $ \sup_{R>0} \|{{\mathcal{M}}}_{{{\mathcal{E}}}_R} \|_{L^p\to L^q} $ are finite.
We let the
{\it sparse exponent set} $\mathrm{Sp}[{{\mathcal{M}}}_{{\mathcal{E}}}] $
consist of all pairs $(\tfrac1{p_1}, \tfrac1{p_2})$ with
$p_2\le p_1'$ for which ${{\mathcal{M}}}_{{\mathcal{E}}}\in \mathrm{Sp}(p_1, p_2) $.
Applying a much more general theorem in \cite{BeltranRoosSeeger-sparse} in our context it then turns out,
that the interiors of these sets are in one-to-one correspondence and we have for all ${{\mathcal{E}}}$
\begin{thm} \label{thm:equiv-impr-sparse}
Let $1<p<q<\infty$. Then
\begin{equation}\label{eq:equiv-impr-sparse}\big(\tfrac 1p,\tfrac 1q\big)\in \mathrm{Int}(\mathscr{L}({{\mathcal{E}}}))
\iff
\big(\tfrac 1p,1-\tfrac 1{q} \big )\in \mathrm{Int}(\mathrm{Sp}[{{\mathcal{M}}}_{{\mathcal{E}}} ])\,.
\end{equation}
\end{thm}

We emphasize that
Theorem \ref{thm:equiv-impr-sparse} implies a result for the interior of the sparse exponent set but not for its boundary.

\subsubsection*{\texorpdfstring{\it Endpoint sparse domination problems}{Endpoint sparse domination}}
The $L^p\to L^q$ estimates for $(p^{-1}, q^{-1})$
are known for the spherical averages
on all open edges for the $L^p\to L^q$ type set
and for the local maximal function ${{\mathcal{M}}}_{[1,2]}$
on all open edges except the vertical edge (where they fail)
(see \cite{StrichartzTAMS1970, Littman}, \cite{Lee03}). However, corresponding nontrivial endpoint sparse domination inequalities for the global variants are open.

\begin{problem}
(i) For the full maximal operator let $P_1,P_2, P_3, P_4$ be as in \eqref{eq:laceysp}. Do we have ${{\mathrm{Sp}}}(p_1,p_2)$ results for $(p_1^{-1}, p_2^{-1})$ on the open boundary edges $(P_1P_4)$ and $(P_4P_3)$ in \eqref{eq:laceysp}?

(ii) For the lacunary maximal operator do we have ${{\mathrm{Sp}}}(p_1,p_2)$ results for $(p_1^{-1}, p_2^{-1})$ on the open edges connecting $(1,0)$ or $(0,1)$ to
$(\frac{d}{d+1}, \frac d{d+1})$?
\end{problem}

\begin{remark} In part (i) sparse domination fails on the open vertical edge $(P_2P_3)$. By \eqref{eq:estmaxform} it would imply an $L^p\to L^{p,\infty}$ inequality for ${{\mathcal{M}}}_{[1,2]}$ for $p=\frac{d}{d-1}$ which fails by the Stein example \eqref{eq:stein-example} (for such an inequality we need to restrict to the smaller space $L^{p,1}$ and even this only works for $d\ge 3$ \cite{STW-jussieu}).
\end{remark}

\begin{remark} For the case $f={{\mathbbm 1}}_F$ and $g={{\mathbbm 1}}_G$ with measurable sets $F,G$ of finite measure, an endpoint estimate for
$\inn{{{\mathcal{M}}} {{\mathbbm 1}}_F}{{{\mathbbm 1}}_G} $ at the corner $P_3=(\frac{d-1}{d},\frac{d-1}{d})$ of the sparse exponent set ${{\mathrm{Sp}}}[{{\mathcal{M}}}]$ was observed by Kesler, Lacey and Mena \cite{KeslerLaceyMena2020}. This does not seem strong enough yet to deduce sparse domination on the $(P_3P_4)$ edge.
\end{remark}

\section{\texorpdfstring{${L}^{{p}}$ improving}{Lp improving}}\label{sec:Lpimpr}
To understand the sparse domination inequalities for the operator ${{\mathcal{M}}}_{{\mathcal{E}}}$ we have to determine the uniform $L^p\to L^q$ bounds for the operators ${{\mathcal{M}}}_{{{\mathcal{E}}}_R}$ for $R>0$.
We then fix a dilation set $E\subset[1,2]$ and consider the problem of $L^p\to L^q$ estimates for the spherical maximal operator ${{\mathcal{M}}}_E$, with the intent to apply some uniform version of these to the sets ${{\mathcal{E}}}_R$.

\subsection{Type sets} \label{sec:typesets}
For $E\subset[1,2]$ define
\[ {{\mathcal{T}}}_E= \big\{ (\tfrac 1p,\tfrac 1q): \,\,\text{ ${{\mathcal{M}}}_E: L^p\to L^q$ is bounded}\big\}\]
as the {\it type set} for ${{\mathcal{M}}}_E$. One would like to determine ${{\mathcal{T}}}_E$ for every given $E$, but given the above mentioned open problems with endpoint estimates even in the case $p=q$ this is difficult. To avoid technicalities we will ignore most endpoint questions and aim to determine the closure of the type set. By convexity, the closure determines the interior.

\begin{problem}\label{typesetproblem}
For given $E\subset [1,2]$ determine the closure $\overline{{{\mathcal{T}}}_E}$ of the type set. \end{problem}

In view of this open problem it makes sense to ask
which subsets of $[0,1]^2$ can arise as $\overline{{{\mathcal{T}}}_E}$ for some $E\subset [1,2]$.
This question was answered by the authors in \cite{RoosSeeger}.
Given parameters $0\le \beta\le \gamma\le 1$ let $\mathcal{Q}(\beta,\gamma)$ be the quadrilateral with corners
\begin{equation} \label{Qbetagamma}
\begin{gathered} Q_1=(0,0), \quad Q_{2,\beta}=(\tfrac{d-1}{d-1+\beta},\tfrac{d-1}{d-1+\beta}),
\\
Q_{3,\beta}=(\tfrac{d-\beta}{d-\beta+1}, \tfrac{1}{d-\beta+1}), \quad Q_{4,\gamma} = (\tfrac{d(d-1)}{d^2+2\gamma-1}, \tfrac{d-1}{d^2+2\gamma-1}).
\end{gathered}
\end{equation}

The quadrangle $\mathcal{Q}(1,1)$ coincides with the closure of the type set for the case $E=[1,2]$ determined by Schlag \cite{Schlag1997} and Schlag--Sogge \cite{SchlagSogge1997}.
For general $E$, it turns out that the closure of the type set is not necessarily a quadrangle.

\begin{thm}[\cite{RoosSeeger}]\label{thm:typesets}
Let $A\subset [0,1]^2$.
There exists $E\subset [1,2]$ such that $\overline{\mathcal{T}_E}=A$ if and only if
$A$ is a closed convex set with
\[ \mathcal{Q}(\beta,\gamma) \subset A\subset \mathcal{Q}(\beta,\beta)\]
for some $0\le \beta\le\gamma\le 1$.
\end{thm}

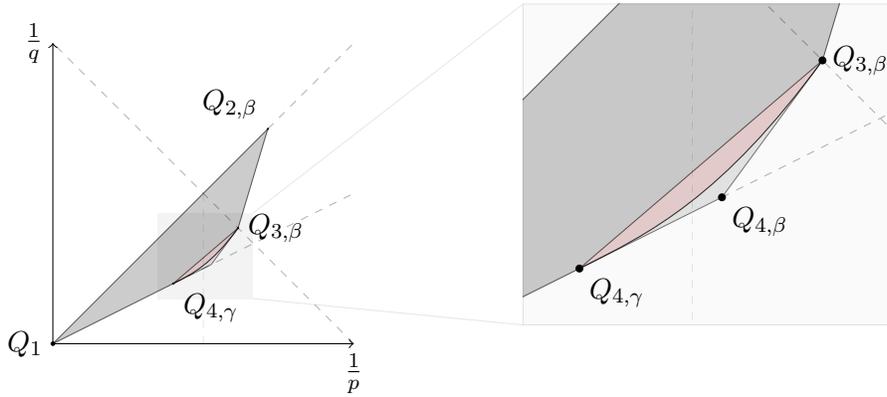
\begin{figure}[ht]
\begin{center}
\begin{tikzpicture}
\def\d{2}
\def\b{.4}
\def\g{1}

\def\clipPad{.05}
\def\clipLx{ \Qfourx{\g} - \clipPad }
\def\clipLy{ \Qfoury{\g} -\clipPad }
\def\clipHx{ \Qthreex{\b} + \clipPad }
\def\clipHy{ \Qthreey{\b} + \clipPad }
\def\clipStyle{solid}
\def\clipOpacity{.1}

\def\clipCx{ (\clipLx+\clipHx)/2 }
\def\clipCy{ (\clipLy+\clipHy)/2 }
\def\cliprad{ .2 }

\begin{scope}[scale=4]
\definecoords

\coordinate (clipLL) at ({\clipLx}, {\clipLy});
\coordinate (clipHL) at ({\clipHx}, {\clipLy});
\coordinate (clipHH) at ({\clipHx}, {\clipHy});
\fill[opacity=.05] (clipLL) rectangle (clipHH);

\drawauxlines{Q4g}{Q2}
\drawQbg
\drawQbb
\drawA
\end{scope}

\begin{scope}[xshift=1cm, yshift=-2cm, scale=15]
\definecoords

\coordinate (clipLLw) at ({\clipLx}, {\clipLy});
\coordinate (clipLHw) at ({\clipLx}, {\clipHy});
\coordinate (clipHHw) at ({\clipHx+.015}, {\clipHy});
\draw[\clipStyle, opacity=\clipOpacity] (clipHL) -- (clipLLw);
\draw[\clipStyle, opacity=\clipOpacity] (clipHH) -- (clipLHw);
\draw[\clipStyle, opacity=\clipOpacity] (clipLLw) rectangle (clipHHw);
\fill[opacity=.02] (clipLLw) rectangle (clipHHw);
\clip (clipLLw) rectangle (clipHHw);

\drawauxlines{Q4b}{Q2}
\drawQbg
\drawQbb
\lblQbb
\drawA

\end{scope}

\end{tikzpicture}
\end{center}
\caption{
Possible type sets.}\label{fig:zoom}
\end{figure}

Figure \ref{fig:zoom} illustrates the situation in the theorem.
Within the critical triangle spanned by the points $Q_{4,\gamma}$, $Q_{4,\beta}$, $Q_{3,\beta}$, the boundary of $\mathcal{T}_E$ may follow an arbitrary convex curve segment.

The parameters $\beta$ and $\gamma$ can be realized as the Minkowski and quasi-Assouad dimensions of $E$, respectively.
To discuss this further we first review notions from fractal geometry.

\subsection{Fractal dimensions}\label{sec:fractaldim}
The
{\em Assouad dimension} of a set ${{\mathcal{E}}}\subset {{\mathbb{R}}}^+$ is defined as
the infimum of all $g>0$ such that there exists $C>0$ such that
\begin{equation}\label{eq:assouaddef} N({{\mathcal{E}}}\cap J, \delta) \le C \big(\tfrac{\delta}{|J|_\times}\big)^{-g}
\end{equation}
holds for all $\delta>0$ and all intervals $J$ with $\delta<|J|_\times<\infty$.
Note that if ${{\mathcal{E}}}=E\subset [1,2]$, then one can restrict to $J\subset [1,2]$ and then $|J|_\times\approx |J|$.

The Assouad dimension can be thought of as the
``largest local Minkowski dimension'' of ${{\mathcal{E}}}$.
We note that $\beta_{{\mathcal{E}}}\le \dim_{{\mathrm{A}}} \!{{\mathcal{E}}}$. For standard self-similar
Cantor sets in $[1,2]$ the two dimensions agree. In contrast, for the sequence example $E= \{1+n^{-a}\,:\,n\in\mathbb{N}\}$ with $a>0$ we have $\beta_{E}= (1+a)^{-1}<1=\dim_{{\mathrm{A}}} E.$

The Assouad dimension itself is not directly relevant for our further discussions.
Instead, a key notion will be the {\it Assouad spectrum}, a refinement of the Assouad dimension introduced by Fraser and Yu \cite{FraserYu2018Adv}.
The idea is to define a continuum of intermediate dimensions interpolating between Minkowski and Assouad dimensions.
For $0\le\theta<1$ one defines $\dim_{{\mathrm{A}},\theta} {{\mathcal{E}}}$ as the infimum of all $g>0$ such that there exists $C>0$ such that
\eqref{eq:assouaddef} holds for all $\delta\in (0,1)$ and all intervals $J$ with diameter $|J|_\times=\delta^\theta$.
The {\it Assouad spectrum} of ${{\mathcal{E}}}$ is the function $\theta\mapsto \dim_{{\mathrm{A}},\theta} {{\mathcal{E}}}$.

A variant is given by $\overline\dim_{{\mathrm{A}},\theta} {{\mathcal{E}}} $ where the condition $|J|_\times=\delta^\theta$ is replaced by $1\ge |J|_\times\ge \delta^\theta$. The function $\theta\mapsto \overline{\dim}_{{\mathrm{A}},\theta}\,{{\mathcal{E}}}$ is the {\it upper Assouad spectrum} of ${{\mathcal{E}}}$
(see Fraser--Hare--Hare--Troscheit--Yu \cite{FraserHareHareTroscheitYu}) and has the benefit of being monotone.

Properties of these intermediate dimensional spectra and many examples can be found in \cite{FraserYu2018Adv, FraserYu2018Ind} and in the monograph by J. Fraser \cite{FraserBook} which also contains many references to further work.
Directly from the definitions one gets for $0\le\theta<1$
\[\beta_{{\mathcal{E}}}= \dim_{{\mathrm{A}},0}{{\mathcal{E}}} \le \dim_{{\mathrm{A}},\theta} {{\mathcal{E}}} \le \overline\dim_{{\mathrm{A}},\theta}\,{{\mathcal{E}}}\le \dim_{{\mathrm{A}}} {{\mathcal{E}}}. \]
The function $\theta\mapsto \dim_{{\mathrm{A}},\theta} {{\mathcal{E}}}$ is continuous on $[0,1)$ and the limit
\[ \dim_{{\mathrm{qA}}} {{\mathcal{E}}}=\lim_{\theta\to 1-} {\dim}_{{\mathrm{A}},\theta}\,{{\mathcal{E}}}\] exists. It is called the {\it quasi-Assouad dimension} (L\"u--Xi \cite{L"uXi2016}). Note that $\dim_{{\mathrm{qA}}} {{\mathcal{E}}}\le \dim_{{\mathrm{A}}}{{\mathcal{E}}}$ and the inequality can be strict. For example, if $E=\{1+ 2^{-\sqrt n}~:~n\in {{\mathbb{N}}}\}$, then $\dim_{{\mathrm{qA}}}E=0$ and $\dim_{{\mathrm{A}}} E=1$.

The sets ${{\mathcal{E}}}$ which have a {\it maximal} Assouad spectrum will be of particular interest to us. For $0\le\beta\le\gamma\le 1$ define
$h_{\beta,{\gamma}}:[0,1)\to [\beta, \gamma] $ by
\begin{equation}\label{hbetagamma}
h_{\beta, {\gamma}}(\theta)= \min(\tfrac{\beta}{1-\theta},\gamma)= \begin{cases}
\frac{\beta}{1-\theta} &\text{ if } 0\le\theta\le 1-\frac \beta\gamma,\\ \gamma &\text{ if }1-\frac \beta\gamma\le \theta<1
\end{cases}
\end{equation}
when $\gamma>0$, and $h_{0,0}\equiv 0$.
Then if $\beta=\beta_{{\mathcal{E}}}$ and $\gamma=\dim_{\mathrm{qA}} {{\mathcal{E}}}$ it follows from the definition of $\gamma$
and the inequality $N({{\mathcal{E}}}\cap J,\delta){\lesssim}_{\varepsilon}\delta^{-\beta-{\varepsilon}}$ that
\begin{equation} \label{eq:lehbg} \dim_{{\mathrm{A}},\theta} {{\mathcal{E}}}\le h_{\beta,\gamma}(\theta), \quad 0\le \theta<1.\end{equation}
We call sets ${{\mathcal{E}}}$
for which we have equality in \eqref{eq:lehbg}{\it quasi-Assouad regular}, or {\it $(\beta,\gamma)$-regular}, if $\beta=\beta_{{\mathcal{E}}}$ and $\gamma=\dim_{\mathrm{qA}} {{\mathcal{E}}}$.
Obvious examples are sets with $\beta_{{{\mathcal{E}}}}=\dim_{{\mathrm{qA}}} {{\mathcal{E}}}$ such as standard Cantor sets. The sets $E=\{1+n^{-a}\,:\,n\in\mathbb{N}\}$ are also quasi-Assouad regular, with $\beta=(1+a)^{-1}$ and $\gamma=1$.
It is natural to ask whether a given function on $[0,1)$ is the Assouad spectrum or the upper Assouad spectrum of some set $E\subset [1,2]$.
Characterizations of both situations were established by Rutar \cite{Rutar24}.

\subsection{Known results}\label{sec:LpLq} We now continue discussing results
on the maximal operators ${{\mathcal{M}}}_E$ from \cite{AHRS, RoosSeeger} for arbitrary dilation sets $E\subset [1,2]$.

\begin{figure}[ht]
\begin{center}
\begin{tikzpicture}[scale=100]
\def\d{4}
\def\b{.3}
\def\g{1}

\def\clipStyle{solid}
\def\clipLineOpacity{1}
\def\clipFillOpacity{0}

\begin{scope}[yscale=1]
\definecoords

\def\ptsize{.015pt}

\def\clipLx{ \Qfourx{\g} - .005 }
\def\clipLy{ \Qfoury{\g} - .006 }
\def\clipHx{ \Qthreex{0} + .01 }
\def\clipHy{ \Qthreey{\b} + .005 }
\def\clipCx{ (\clipLx+\clipHx)/2 }
\def\clipCy{ (\clipLy+\clipHy)/2 }
\def\cliprad{ .2 }

\coordinate (clipLL) at ({\clipLx}, {\clipLy});
\coordinate (clipHH) at ({\clipHx}, {\clipHy});
\draw[\clipStyle, opacity=\clipLineOpacity] (clipLL) rectangle (clipHH);
\fill[opacity=\clipFillOpacity] (clipLL) rectangle (clipHH);
\clip (clipLL) rectangle (clipHH);

\drawauxlines{Q4b}{Q2}

\def\QbgDrawCritSeg{1}
\def\QbgCritSegStyle{solid}
\def\QbgCritSegOpacity{.2}
\def\AFillColor{black}
\def\QbbFillOpacity{0}
\def\QbgFillOpacity{.1}
\def\AFillOpacity{.1}

\drawQbg
\drawQbb
\lblQbb
\drawA
\lblQthreezero

\def\tangentOpacity{.8}
\newcommand{\drawATangent}[3]{
\path (Q4g) .. controls (C1) and (C2) .. (Q3) node [sloped,pos=#1,minimum width=10cm](tang) {};
\coordinate (tmpA) at (intersection of tang.west--tang.east and Q4g--Q30);
\coordinate (tmpB) at ($ (intersection of tang.west--tang.east and Q3--Q30) + (.011pt,.011pt) $);
\draw[opacity=\tangentOpacity] (tmpA) -- (tmpB);
\fill (tmpA) node [below right] {#2} circle [radius=\ptsize];
\fill (tmpB) node [right] {#3} circle [radius=\ptsize];
}

\def\ptsize{.01pt}
\foreach \x[evaluate=\x as \p using \x*0.1] in {2,...,9} {
\ifnum\x=5
\drawATangent{\p}{$Q_{4,\gamma_n}$}{$Q_{3,\beta_n}$}
\else
\drawATangent{\p}{}{}
\fi
}

\end{scope}
\end{tikzpicture}
\end{center}
\caption{Building a type set.}

\label{fig:tangents}

\end{figure}
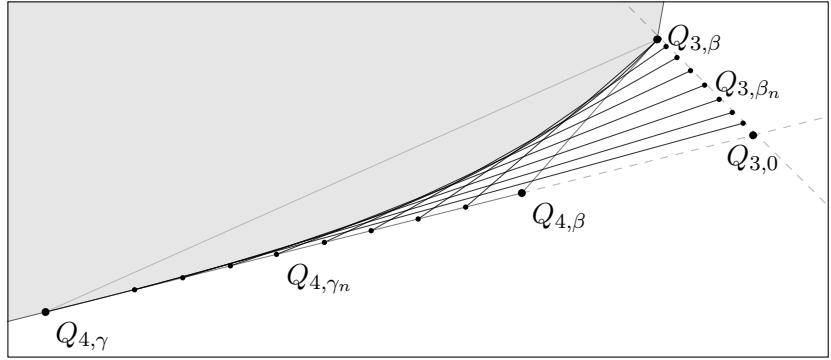

\begin{thm}[\cite{RoosSeeger}]\label{thm:betagammath} Let $d\ge 2$.
Then the following hold.

(i) For all $E\subset[1,2]$ with
$\beta=\beta_E$, $\gamma=\dim_{{\mathrm{qA}}}\,E$ we have
\[ \mathcal{Q}(\beta,\gamma)\subset\overline{\mathcal{T}_E} \subset {{\mathcal{Q}}}(\beta,\beta).\]

(ii) If $E$ is quasi-Assouad regular, then \[\overline{\mathcal{T}_E}=\mathcal{Q}(\beta,\gamma). \]

(iii) Let $E=\bigcup_{k=1}^N E_k$, where $E_k$ is a $(\beta_k, {\gamma}_k)$-regular set, for $k=1,\dots, N$. Then
\[ \overline{{{\mathcal{T}}}_E} = \bigcap_{k=1}^N {{\mathcal{Q}}}(\beta_k, \gamma_k).\]
\end{thm}
Part (iii) is an immediate consequence of part (ii), and we observe that for a finite union of quasi-Assouad regular sets the closure of the type set is a polygonal region.
The construction of $E$ in the type set characterization of Theorem \ref{thm:typesets} relies on part (iii):
one carefully builds $E$ as a countable union of suitable quasi-Assouad regular sets with dimensional parameters $(\beta_k,\gamma_k)$ so that the sets
$\bigcap_{k=1}^N {{\mathcal{Q}}}(\beta_k,\gamma_k)$ approximate the set $A$ (see Figure \ref{fig:tangents}).

\begin{proof}[Outline of proof for the inclusion ${{\mathcal{Q}}}(\beta,\gamma)\subset \overline{{{\mathcal{T}}}_E}$] One analyzes the $L^p\to L^q$ bounds for the maximal operator $M_E^j$ in \eqref{eq:Mj-def}.
A straightforward $O(2^j)$ bound for the convolution kernel $\sigma_j$ yields
\begin{equation}\label{eq:1infty} \|M^j_E\|_{L^1 \to L^\infty} {\lesssim} 2^j.\end{equation}
Interpolating this with the $L^2$ bound $\|M^j_E\|_{L^2\to L^2} {\lesssim}_{\varepsilon} 2^{-j (\frac{d-1-\beta-{\varepsilon}}{2}) }$ (cf. \eqref{eq:p'bd}) yields the asserted behavior for $(\frac 1p, \frac 1q)$ near $Q_{3,\beta}$.

The estimates for $(\frac 1p,\frac 1q)$ close to $Q_{4,{\gamma}}$ are more interesting. The objective is to develop versions of the space-time estimates in \cite{SchlagSogge1997} that apply in a fractal situation.
Here we distinguish the cases
\begin{itemize} \item[(i)] $d\ge 3$,
or $d=2$
and ${\gamma}<1/2$, \item[(ii)] $d=2$, ${\gamma}\ge 1/2$.
\end{itemize} In the first case one interpolates the $L^1\to L^\infty$ bound \eqref{eq:1infty} with
an $L^2 \to L^{q}$ bound
\[\|M^j_E\|_{L^2 \to L^{q_{\gamma}} } {\lesssim} C_{\varepsilon} 2^{j ({\varepsilon}- \frac{(d-1)^2-2\gamma}{2(d-1)+2{\gamma}})} \quad \text{ for } q_{\gamma}= \tfrac{2(d-1+2{\gamma})}{d-1}\]
which itself relies on a fractal version of a Strichartz type estimate \cite{AHRS}
\[ \Big(\sum_{t\in E(j)} \|f*\sigma^j_t\|_ {q_{\gamma}}^{q_{\gamma}} \Big)^{1/q_{\gamma}} {\lesssim}_{\varepsilon}
2^{j ({\varepsilon}- \frac{(d-1)^2-2\gamma}{2(d-1)+2{\gamma}})}\|f\|_2.
\] Here $E(j)$ is a maximal $2^{-j}$-separated subset of $E$.
For the interpolation to work the $L^2\to L^q$ operator norm needs to exponentially decay as $j\to \infty$ and for this we need to have $(d-1)^2-2{\gamma}>0$, corresponding to case (i).

In case (ii) the proof is substantially harder, and we refer to \cite{RoosSeeger}. One uses bilinear operators and the argument relies on some fractal version of a theorem by Klainerman and Machedon \cite{KlainermanMachedon}.
\end{proof}
Various endpoint results for $L^p\to L^q$ boundedness on the edges of ${{\mathcal{Q}}}(\beta,\gamma)$, with
appropriate assumptions on $N(E\cap J,\delta)$, $N(E,\delta)$,
have been obtained in \cite{AHRS}, \cite{RoosSeeger}, \cite{LeeRoncalZhangZhao}. These results tend to be sharp for quasi-Assouad regular $E$, but not in general.

One can combine Theorem \ref{thm:equiv-impr-sparse} with an $R$-uniform version of Theorem \ref{thm:betagammath} applied to the operators ${{\mathcal{M}}}_{{{\mathcal{E}}}_R}$
to obtain a sparse domination result.
The required $R$-uniform version follows from an examination of the proof of part (i) of Theorem \ref{thm:betagammath}.
Assume that ${{\mathcal{E}}}\subset {{\mathbb{R}}}^+$, $\beta=\beta_{{\mathcal{E}}}$ and $\gamma=\dim_{q{\mathrm{A}}} {{\mathcal{E}}}$.
Let $(\tfrac 1{p_1},\tfrac 1{p_2} )$ belong to the interior of the quadrilateral ${{\mathcal{P}}}_{\mathrm{sp}}(\beta,\gamma)$ with corners
\begin{equation} \label{Pbetagamma}
\begin{gathered} P_1=(0,1), \quad P_2(\beta)=(\tfrac{d-1}{d-1+\beta},\tfrac{\beta}{d-1+\beta}),
\\
P_3(\beta)=(\tfrac{d-\beta}{d-\beta+1}, \tfrac{d-\beta}{d-\beta+1}), \quad P_4(\gamma) = (\tfrac{d(d-1)}{d^2+2\gamma-1}, \tfrac{d^2-d+2{\gamma}}{d^2+2\gamma-1}).
\end{gathered}
\end{equation}
Then ${{\mathcal{M}}}_{{\mathcal{E}}}\in {{\mathrm{Sp}}}(p_1,p_2)$.
The result is sharp up to the boundary of ${{\mathcal{P}}}_{\mathrm{sp}}(\beta,\gamma)$ if ${{\mathcal{E}}}=\cup_k2^kE$ with $E$ a $(\beta,\gamma)$-regular set. One can use uniform versions of part (iii) of Theorem \ref{thm:betagammath} to formulate other essentially sharp sparse domination results.

\subsection{Necessary conditions}
We now summarize necessary conditions which will suggest a conjecture. Let $(\frac 1p,\frac 1q)\in \overline{{\mathcal{T}}}_E$. From general considerations for convolution operators we must have $p\le q$, i.e. $(\frac 1p,\frac 1q)$ is on or below the line connecting $Q_1$ and $Q_{2,\beta}$. Since for $t\in E$ we have ${{\mathcal{M}}}_E f\ge {{\mathcal{M}}}_{\{t\}}f$ we may test on a characteristic function of a thin $\delta$-annulus with radius $t$ and find that
$\frac{d}{q}\ge \frac 1p$; i.e. $(\frac 1p,\frac 1q)$ must lie on or above the line through $Q_1$ and $Q_{4,\gamma}$.
By testing the operator on characteristic functions of balls of radius $\delta$ (arguing as before for the case $p=q$) we find that
\[\delta^{d-1+\frac 1q} N(E,\delta)^{\frac 1q} {\lesssim} \|{{\mathcal{M}}}_E\|_{p\to q}\, \delta^{\frac dp} \] and
since for every ${\varepsilon}>0$ there are arbitrarily small $\delta$ such that $N(E,\delta)\gtrsim \delta^{-\beta+{\varepsilon}}$, the point $(\frac 1p,\frac 1q)$ must lie on or to the left of the line through $Q_{2,\beta}$ and $ Q_{3,\beta}$.

A more interesting observation is that for
every $\delta>0$, and for all intervals $J$ with $\delta\le |J|\le 1$ we must have
\begin{equation}\label{eq:newneccond} N(E\cap J,\delta)^{\frac 1q} {\lesssim} \|{{\mathcal{M}}}_E\|_{p\to q} \, \delta^{\frac 1p-\frac dq} \Big(\frac{\delta}{|J|}\Big)^{\frac{d-1}{2}(\frac 1p+\frac 1q-1)} .
\end{equation}
To check \eqref{eq:newneccond}
one tests the operator on a characteristic function of a $\delta$-neighborhood of a spherical cap of diameter $\sqrt \rho$, with $\rho={\delta}/{|J|}$ (so that $\sqrt \rho$ is intermediate between $\sqrt \delta$, corresponding to the usual Knapp example, and $1$, which is essentially the thin annulus example above). If $J=[a,b]$, one computes $ {{\mathcal{M}}}_E f_{\delta,\rho} {\gtrsim} \rho^{\frac{d-1}{2}}$ for $|x'|\le c \frac{\delta}{\sqrt\rho}$, $t\in E\cap J$, $|x_d+t-a| \le c\delta$, which then leads to \eqref{eq:newneccond} (see \cite[\S5]{RoosSeeger}).

Using intervals $J$ of length $\approx 1$ one checks that $(\frac 1p,\frac 1q)\in {{\mathcal{T}}}_E$ cannot lie below the line through $Q_{4,\beta}$ and $ Q_{3,\beta}$, and a combination of the above leads to the necessary condition $\overline {{\mathcal{T}}}_E\subset {{\mathcal{Q}}}(\beta,\beta)$ in part (i) of the theorem.

If $E$ is $(\beta,\gamma)$-regular we choose suitable $J$ with $| J|=\delta^\theta$ where $\theta=1-\frac{\beta}{\gamma}$ so that essentially
$N(E\cap J ,\delta){\gtrsim} \delta^{-(1-\theta)\gamma+{\epsilon}}{\gtrsim} \delta^{-\beta+{\epsilon}'}$
and a computation shows that $(\frac 1p,\frac 1q)$ cannot lie below the line
through $Q_{4,\gamma}$ and $ Q_{3,\beta}$, i.e. the necessary condition from \cite{AHRS} in this case. Thus for $(\beta,{\gamma})$-regular dilation sets the necessary condition is sharpened to $\overline{{{\mathcal{T}}}_E}\subset {{\mathcal{Q}}}(\beta,\gamma)$ in part (ii) of the theorem.
Finally if $E=\bigcup_{k=1}^N E_k$ we have ${{\mathcal{M}}}_E f\ge{{\mathcal{M}}}_{E_k} f$ which gives the inclusion $\overline{{{\mathcal{T}}}_E}\subset {{\mathcal{Q}}}(\beta_k,{\gamma}_k) $ for $k=1,\dots, N$ i.e. the necessary condition for part (iii) of the theorem.

\subsection{A conjecture}\label{sec:aconjecture}
For general $E\subset[1,2]$ with $\dim_{{\mathrm{M}}}E=\beta$, $\dim_{{\mathrm{qA}}}E=\gamma$ Theorem \ref{thm:betagammath}
leaves open what happens in the triangle ${{\mathcal{R}}}(\beta,\gamma)={{\mathcal{Q}}}(\beta,\beta)\setminus {{\mathcal{Q}}}(\beta,\gamma)$. A reasonable conjecture for the part of $\overline{{{\mathcal{T}}}_E}$ in this triangle is that the necessary condition \eqref{eq:newneccond} for membership in ${{\mathcal{T}}}_E$ is also sufficient for $(\frac 1p,\frac 1q)$ to belong to $\overline{{{\mathcal{T}}}_E}$. In other words, we conjecture that for $(\frac 1p,\frac 1q)\in {{\mathcal{R}}}(\beta,\gamma)$ the pair $(\frac 1p,\frac 1q)$ belongs to $\overline{{{\mathcal{T}}}_E}$ if
\begin{equation}\label{eq:genLpLQconj}\sup_{\delta\le |J|\le 1} N(E\cap J,\delta) |J|^{-\frac{q(d-1)}{2}(1-\frac 1p-\frac 1q)} {\lesssim} \delta^{-\big(\frac{(d-1)q}2 - \frac{(d+1)q}{2} (\frac 1p-\frac 1q) \big)} \end{equation}
where the supremum is taken over all intervals $J$ with length between $\delta$ and $1$.
To formulate this as a conjectural equivalence for $\overline{{{\mathcal{T}}}_E}$
we use another dimensional quantity which was introduced in \cite{BeltranRoosSeeger-radial,BeltranRoosRutarSeeger}.
For an arbitrary set ${{\mathcal{E}}}\subset {{\mathbb{R}}}^+$ the {\it Legendre--Assouad function} is defined by
\begin{equation}\label{eq:LA} \nu_{{\mathcal{E}}}^\sharp(\alpha)= \limsup_{\delta\to 0} \frac{\log\big( \sup_{\delta\le |J|_\times\le 1} |J|_\times^{-\alpha} N({{\mathcal{E}}}\cap J, \delta)\big)}{ \log(\delta^{-1})}
\end{equation}
for every $\alpha\in\mathbb{R}$. If ${{\mathcal{E}}}=E\subset [1,2]$, then we may restrict to $J\subset [1,2]$ and note $|J|_\times\approx |J|$.

Recall from Theorem \ref{thm:betagammath} that $\overline{\mathcal{T}_E} \subset {{\mathcal{Q}}}(\beta,\beta)$.
\begin{problem} \label{LAconj}
Let $E\subset[1,2]$ and $\beta=\beta_E$. Let $(\tfrac 1p,\tfrac 1q)\in {{\mathcal{Q}}}(\beta,\beta)$.
Prove that $(\frac 1p,\frac 1q)\in \overline{{{\mathcal{T}}}_E}$ if and only if
\begin{equation} \label{eqLnex-LA} \tfrac 1q \nu_E^\sharp\big( \tfrac{q(d-1)}{2}(1-\tfrac 1p-\tfrac 1q)\big) \le \tfrac{d-1}{2} - \tfrac{d+1}{2}(\tfrac 1p-\tfrac 1q)\,.
\end{equation}
\end{problem}

As proved in the recent paper \cite{BeltranRoosRutarSeeger} by Beltran, Rutar and the authors, the function $\nu_{{\mathcal{E}}}^\sharp$
in \eqref{eq:LA}
turns out to be the {\it Legendre transform} of the function \[\nu_{{\mathcal{E}}}(\theta)=-(1-\theta)\mathrm{dim}_{\mathrm{A},\theta}\,{{\mathcal{E}}}, \] that is,
\begin{equation} \label{eq:ass-leg}\nu_{{\mathcal{E}}}^\sharp (\alpha)= \sup_{0\le \theta\le 1} \big ( \alpha\theta-\nu_{{\mathcal{E}}}(\theta)\big).
\end{equation}
The following lemma lists some elementary properties (see \cite{BeltranRoosSeeger-radial,BeltranRoosRutarSeeger, FraccaroliRoosSeeger}).
\begin{lemma} \label{lem:nusharp-prop} Let $E\subset [1,2]$ with $\dim_{\mathrm{M}} E=\beta$ and $\dim_{{\mathrm{qA}}}E=\gamma$. Then the following hold.

(i) $\nu_E^\sharp$ is convex and nondecreasing on ${{\mathbb{R}}}$.

(ii) There exists an $\alpha_*\in [0,\gamma]$ such that
$\nu_E^\sharp(\alpha)= \beta$ for $\alpha\le \alpha_*$ and such that $\nu_E^\sharp$ is strictly increasing on $(\alpha_*,\infty)$.

(iii) $\nu_E^\sharp(\alpha)=\alpha$ for $\alpha\ge \gamma$.
\end{lemma}

Part (i) follows from properties of the Legendre transform. Note that the conjectured range \eqref{eqLnex-LA} in Problem \ref{LAconj}
is consistent with the known inclusion ${{\mathcal{Q}}}(\beta,\gamma)\subset\overline{{{\mathcal{T}}}_E}$, by part (iii) of the lemma. We also recover
the known results for quasi-Assouad-regular sets or finite unions of those:

\begin{examples}

(i) Let $E\subset [1,2]$ be $(\beta,\gamma)$-Assouad regular, i.e. with maximal Assouad spectrum $\dim_{\mathrm A,\theta}E= h_{\beta,\gamma}(\theta)$ as in \eqref{hbetagamma}.
Then $\nu_E^\sharp $ is given by
\[\nu_{{\mathcal{E}}}^\sharp(\alpha)= \begin{cases}
\beta &\text{ if } \alpha\le 0,
\\
(1-\tfrac \beta\gamma)\alpha+\beta &\text{ if $0\le \alpha\le \gamma $, }
\\
\alpha &\text{ if } \alpha>\gamma.
\end{cases}
\]

(ii) For $\ell=1,\dots, N$ let $E^\ell\subset [1,2]$, $\ell=1,\dots, N$ be $(\beta_\ell,\gamma_\ell)$-regular, and let $E=\cup_{\ell=1}^N E^\ell$.
Then
\[ \nu_E^\sharp(\alpha)=\max_{\ell=1,\dots, N} \nu_{E^\ell}^\sharp(\alpha).\]
\end{examples}

It is natural to ask which functions on $[0,\infty)$ can arise as $\nu^\sharp_E|_{[0,\infty)}$ for suitable $E\subset [1,2]$.
\begin{proposition}[\cite{BeltranRoosRutarSeeger}]\label{prop:whichhassleg}
Consider a function $\tau:[0,\infty)\to[0,\infty)$. There exists a set $E$ such that $\tau(\alpha)=\nu_E^\sharp(\alpha)$ for all $\alpha\ge 0$ if and only if $\tau$ is increasing, convex, and satisfies $\tau(\alpha)=\alpha$ for $\alpha\ge 1$.
\end{proposition}

\begin{remark}
The identity \eqref{eq:ass-leg} remains true if one replaces the Assouad spectrum by the upper Assouad spectrum in the definition of $\nu_E$. This is because the Legendre transform only depends on the convex hull of the input function and both these functions have the same convex hull.
The Legendre transform also appears in connection with the Assouad spectrum in
a different context in the fractal geometry paper \cite{BFKR} by Banaji, Fraser, Kolossv\'ary and Rutar.
\end{remark}

\begin{remark}
One can also ask about the $L^p$ improving inequalities for ${{\mathcal{M}}}_E$ when acting on radial functions. This was considered in our paper with David Beltran \cite{BeltranRoosSeeger-radial} and also by Shuijiang Zhao \cite{ZhaoSh24}.
In dimension $d\ge 3$ one gets a very simple and clean characterizations
of $L^p_{{\text{\rm rad}}} \to L^q$ boundedness; however in view of the nonradial example discussed above these estimates are not useful anymore to predict sharp outcomes for the general case. This comment also applies in two dimensions, but the analysis in the radial case now involves a weakly singular kernel, which makes the 2D case
more interesting.
In \cite{BeltranRoosSeeger-radial}
$\nu_E^\sharp$ was used to characterize the closure of the type set
\[{{\mathcal{T}}}^{\mathrm{rad,2}}_E= \big\{ (\tfrac 1p,\tfrac 1q): \text{ ${{\mathcal{M}}}_E$ maps $L^p_{{\text{\rm rad}}}({{\mathbb{R}}}^2)$ to $L^q({{\mathbb{R}}}^2)$} \big\}. \] Specifically, if $\beta=\dim_{{\mathrm{M}}}E$ and $\Delta_\beta$ is the triangle with corners
$(0,0)$, $(\tfrac{1}{1+\beta}, \tfrac 1{1+\beta})$ and $(\tfrac{2}{3+\beta}, \tfrac{1}{3+\beta})$ then in two dimensions
\[ \overline{\mathcal{T}_E^{\mathrm{rad,2}}} = \Delta_\beta \, \cap \, \big\{ (\tfrac1p,\tfrac1q): \tfrac1q \nu_E^\sharp(\tfrac{q}2-1)+\tfrac1p - \tfrac1q \le \tfrac 12 \big\}. \]
\end{remark}

\begin{remark}
When considering $L^p\to L^q$ estimates for $p<q$ the restriction to dilation sets to a compact subset of ${{\mathbb{R}}}^+$ is mandated by scaling considerations. In order to extend to global dilation sets on has to necessarily include a scaling factor in the definition of the maximal operator, and work with the ``fractional'' modification
\[{{\mathcal{M}}}_{{{\mathcal{E}}},\nu} f(x) = \sup_{t\in {{\mathcal{E}}}} t^{\nu} |A_tf(x)|,\text{ with } \nu=\tfrac dp-\tfrac dq. \]
In particular, interesting cases where we know $L^p\to L^q$ boundedness for ${{\mathcal{M}}}_E$ with $E\subset [1,2]$ lead to endpoint problems for the operator ${{\mathcal{M}}}_{{{\mathcal{E}}},\nu}$ with ${{\mathcal{E}}}=\cup_{k\in {{\mathbb{Z}}}} 2^k E$ and $\nu=d/p-d/q$.
Here we refer to work in \cite{Oberlin89}.
Recently, restricted weak type $(p,q)$ endpoint estimates for the full modified spherical maximal functions have been obtained in work by Basak, Choudhary and Spector \cite{BasakChoudharySpector}. These go along with the open question whether such restricted weak type estimates can be upgraded to strong type $(p,q)$ estimates. Interestingly, in some problems where one works with {\it nonisotropic} dilations such global extensions are easier to obtain via bootstrapping arguments of Nagel--Stein--Wainger type \cite{NagelSteinWainger-lac}, under an assumptions of disjointness of the wavefront set of the underlying measure and the eigenspaces of the dilation operators. For suitable examples this had been first shown in unpublished work by M. Christ \cite{Christ-frac}, with further variants and discussion in \cite[\S4]{GreenleafSeegerWainger}, \cite{SeegerWainger-frac}.
However it is not clear how to apply these ideas to endpoint bounds in the presence of standard isotropic dilations.
\end{remark}

\section{On power weight inequalities}\label{sec:power}
Here we discuss weighted norm inequalities with power weights $w_\alpha(x)=|x|^{\alpha}$ (already mentioned in \S\ref{sec:sparse}). The norm in $L^p(w_\alpha)$ is given by \[\|g\|_{L^p(w_\alpha)}=\Big(\int|g(x)|^p |x|^\alpha dx \Big)^{1/p} .\] We shall see that
the Legendre--Assouad function of dilation sets plays again a crucial role in understanding the $L^p(w_\alpha)$ boundedness of the spherical maximal operators.

Optimal inequalities for the lacunary maximal spherical maximal function and almost optimal inequalities for the full maximal function ${{\mathcal{M}}}_{\mathrm{full}}$ (with ${{\mathcal{E}}}={{\mathbb{R}}}^+$) are due to Duoandikoetxea and Vega \cite{DuoandikoetxeaVega} who showed that is bounded on $L^p(w_\alpha)$ if
$1-d<\alpha<(d-1)p-d$ and $p>1+\frac{1}{d-1}$. Moreover the conditions $\alpha<(d-1) p-d$ and $\alpha\ge 1-d$ are necessary.

Juyoung Lee \cite{LeeJ-weighted} established the endpoint inequality with $\alpha=1-d$ if $d\ge 3$, $p\ge 2$ and $d=2$ and $p>2$, thus settling the endpoint problem in two dimensions. His argument relies crucially on an application of the inequality
$\sum_j\|P_j f\|_p^p {\lesssim} \|f\|_p^p $ for dyadic frequency decompositions $P_j$ which fails for $p<2$. Consequently the endpoint result remains open for $p<2$:

\begin{problem} Let $d\ge 3$ and $\frac{d}{d-1}<p<2$. Does ${{\mathcal{M}}}_{\mathrm{full}}$ map $L^p(w_{1-d})$ to $L^p(w_{1-d})$?
\end{problem}

We turn to the problem of power weight inequalities for the restricted maximal operators ${{\mathcal{M}}}_{{\mathcal{E}}}$ which was first studied by Duoandikoetxea and Seijo \cite{DuoandikoetxeaSeijo}. Counterexamples from the unweighted theory show that such inequalities fail for $p<p_\beta=1+\frac{\beta}{d-1}$, where $\beta=\beta_{{\mathcal{E}}}$.
Let
\[ {{\mathcal{W}}}_{{\mathcal{E}}}=\big \{ (\tfrac 1p, \tfrac \alpha p)\in [0,1]\times \mathbb{R} \;: \; {{\mathcal{M}}}_{{\mathcal{E}}} \text{ is bounded on } L^p(w_\alpha) \big \}.\]
By analytic interpolation, ${{\mathcal{W}}}_{{\mathcal{E}}}$ is a convex set. Examples from the unweighted theory show that $p\ge 1+\frac{\beta}{d-1}$ is necessary for boundedness on $L^p(w_\alpha)$ to hold.

Duoandikoetxea and Seijo \cite{DuoandikoetxeaSeijo} obtained an essentially sharp result on ${{\mathcal{W}}}_{{\mathcal{E}}}$ for the special case
that $\dim_{A}\,{{\mathcal{E}}}=\beta_{{\mathcal{E}}}$ coincide. Moreover
they exhibited examples which show that in general it does not suffice to just consider the Minkowski dimension for the formulation of essentially sharp results. However for most dilation sets in $[1,2]$ the problem of determining the closure of ${{\mathcal{W}}}_{{\mathcal{E}}}$ was left open. In our recent work with Marco Fraccaroli \cite{FraccaroliRoosSeeger} we determined the closure of ${{\mathcal{W}}}_{{\mathcal{E}}}$ for all dilation sets
by using the Legendre--Assouad function.

\begin{figure}[ht]
\begin{tikzpicture}
\def\d{3}
\def\b{.5}

\def\clipPad{.05}
\def\clipLx{ \Aonex - \clipPad }
\def\clipLy{ \Pfoury{\b} -\clipPad }
\def\clipHx{ \Pfourx{\b} + \clipPad }
\def\clipHy{ \Aoney + \clipPad }
\def\clipStyle{solid}
\def\clipOpacity{.1}

\begin{scope}[scale=3]
\definecoordsmod

\coordinate (clipLL) at ({\clipLx}, {\clipLy});
\coordinate (clipHL) at ({\clipHx}, {\clipLy});
\coordinate (clipHH) at ({\clipHx}, {\clipHy});
\fill[opacity=.05] (clipLL) rectangle (clipHH);

\drawauxlinesmod{P3}
\drawWbg
\drawAmod
\end{scope}

\begin{scope}[xshift=-3cm, yshift=11cm, scale=15]
\definecoordsmod

\coordinate (clipLLw) at ({\clipLx}, {\clipLy});
\coordinate (clipLHw) at ({\clipLx}, {\clipHy});
\coordinate (clipHHw) at ({\clipHx+.015}, {\clipHy});
\draw[\clipStyle, opacity=\clipOpacity] (clipHL) -- (clipLLw);
\draw[\clipStyle, opacity=\clipOpacity] (clipHH) -- (clipLHw);
\draw[\clipStyle, opacity=\clipOpacity] (clipLLw) rectangle (clipHHw);
\fill[opacity=.02] (clipLLw) rectangle (clipHHw);
\clip (clipLLw) rectangle (clipHHw);

\drawauxlinesmod{P3}
\drawWbg
\drawAmod

\end{scope}

\end{tikzpicture}

\caption{
Vertically compressed shape of the closure of a type set ${{\mathcal{W}}}_{{\mathcal{E}}}$ in \eqref{eq:WEdescription} (shaded), with parameters $\gamma=1$,
$x_1=\frac{d-1}{d}$, $x_\beta=\frac{d-1}{d-1+\beta} $.
}

\label{fig:WE}

\end{figure}

To describe this result for ${{\mathcal{E}}}\subset {{\mathbb{R}}}^+$ let
$(\nu^\sharp)^\dagger$ denote the inverse of
the restriction of $\nu_E^\sharp$ to $(\alpha_*,\infty)$ ({\it cf.} Lemma \ref{lem:nusharp-prop}).
It is shown in \cite{FraccaroliRoosSeeger} that
\begin{equation} \label{eq:WEdescription}
\overline{{{\mathcal{W}}}_{{\mathcal{E}}}} = \big \{ \big(\tfrac 1p,\tfrac\alpha p\big )\in [0,1]\times \mathbb{R}\,:\, p\ge 1+\tfrac{\beta}{d-1} ,\,\, L(p)\le \alpha\le U(p) \big\}
\end{equation} where for $p\ge 1+\tfrac{\beta}{d-1}$ the functions $U$ and $L$ are given by
\begin{align*} U(p)&= (d-1)(p-1)-\beta,
\\L(p) &=(d-1)(p-2)- (\nu^\sharp_{{\mathcal{E}}})^{\dagger} \big((d-1)(p-1) \big).
\end{align*}
See Figure \ref{fig:WE} for a vertically compressed copy of a typical $\overline{{{\mathcal{W}}}_{{\mathcal{E}}}}$.
It would also be interesting to address, for specific classes of dilation sets, open questions on endpoint $L^p(w_\alpha)$ estimates for ${{\mathcal{M}}}_{{\mathcal{E}}}$, for both the local and the global versions.

\section{Maximal functions: Further directions}\label{sec:further}

\subsection{Multiparameter spherical maximal functions}
\label{sec:maxfunctionapproachset}
In the pointwise convergence problem for families of spherical means
one can replace the vertical approach set of radii by a more general set $V$ in ${{\mathbb{R}}}^{d}\times {{\mathbb{R}}}^+$ with accumulation point at $(0,0)$. One then asks under which conditions on $V$ we have for all $f\in L^p({{\mathbb{R}}}^d)$
\[\lim_{\substack{(y,t)\to(x,0) \\ (y,t)\in (x,0)+V}} A_t f(y)=f(x) \] for almost every $x\in {{\mathbb{R}}}^d$.
This leads to consideration of
the maximal operator ${{\mathcal{M}}}_V^*$ given by
\[{{\mathcal{M}}}^*_V f(x)= \sup_{(u,t)\in V} \Big|\int f(x+u-t y')\, d\sigma(y') \Big|\]
\begin{problem} Given $V\subset {{\mathbb{R}}}^d\times {{\mathbb{R}}}^+$ find satisfactory $L^p$, $L^p$-improving, sparse domination and weighted $L^p$ bounds for the maximal operators ${{\mathcal{M}}}^*_V$.
\end{problem}
This covers the situation in Wolff's work on the spherical Nikodym maximal function \cite{Wolff-Kakeya-circles} (where $t$ is kept constant). Many other results in this direction can be found in the recent paper by Chang, Dosidis and Kim \cite{ChangDosidisKim}.

Another natural type of multi-parameter maximal function is obtained by replacing the standard dilations by a multiparameter dilation; this leads to ellipsoidal maximal functions, see \cite{Cho-multiparametersph}, \cite{Erdogan2003}.
Significant progress in this direction has been recently achieved based on the phenomenon of {\it multiparameter cinematic curvature}, in works by Lee--Lee--Oh \cite{LeeLeeOh-JFA2025,LeeLeeOh-PAMS}, Chen--Guo--Yang \cite{ChenGuoYang2023}, Zahl \cite{Zahl2012, Zahl2026}, and Hickman--Zahl \cite{HickmanZahl2025}. See also Schippa \cite{Schippa-sqfct2025} for related work on multiparameter local smoothing bounds, via an approach using multiparameter square-functions.

\subsection{Helical averages}\label{sec:helical}
The problems in this survey can all be formulated for other families of averages, even in non-convolution situations, with outcomes depending on the geometry of the wavefront sets of the associated distribution kernels. The model case of the {\it helical averages} in three dimensions
\begin{equation} \label{eq:helicalave} {A^{\mathrm{hel}}}_t f(x)= \frac{1}{2\pi} \int_0^{2\pi} f(x_1-t \cos \alpha, x_2-t\sin\alpha, x_3-t\alpha)\,d\alpha
\end{equation} is of particular interest, as is the {\it helical maximal function} ${M^{\mathrm{hel}}} f(x)=\sup_{t>0}| {A^{\mathrm{hel}}}_t f(x)|$. The optimal $L^p$ improving estimates for ${A^{\mathrm{hel}}}_t$ were proved by Oberlin \cite{Oberlin-helix}, in particular the endpoint $L^{3/2}\to L^{2}$ and $L^2\to L^3$ estimates. This was generalized to a class of Fourier integral operators with one-sided fold singularities in \cite{GreenleafSeeger1994} satisfying a suitable curvature condition in the fibers of the fold surface (reminiscent of the Sogge cinematic curvature condition).
The first sharp $L^p$-Sobolev bounds for the helical averages gaining $1/p$ derivatives
in \cite{PramanikSeeger} used Wolff's decoupling result \cite{Wolff2000} for the light cone as a black box. A combination of \cite{PramanikSeeger} with the optimal Bourgain-Demeter decoupling \cite{BourgainDemeter2015} yields that ${A^{\mathrm{hel}}}_t$ maps $L^p$ to $L^p_{1/p}$ for $p>4$, see also \cite{PramanikSeeger2021} for a variable coefficient generalization. The range $p>4$
is optimal as shown by an example of Oberlin--Smith \cite{OberlinSmith}). Alternatively, the sharpness can also be understood in terms of a Wolff example in decoupling (\cite{Wolff2000}, \cite[\S3]{BGHS-Adv}).
Decoupling arguments were also used
to establish $L^p$-boundedness of the maximal operator ${M^{\mathrm{hel}}}$ in the range $p>4$ (\cite{PramanikSeeger, BourgainDemeter2015}) which however is not optimal. The sharp result for ${M^{\mathrm{hel}}}$
was obtained by Beltran, Guo, Hickman and the second author \cite{BeltranGuoHickmanSeeger}, and independently by Ko, Lee and Oh \cite{KoLeeOh2022}:
\begin{thm} \label{thm:helicalmax} ${M^{\mathrm{hel}}}$ is bounded on $L^p({{\mathbb{R}}}^3)$ if and only if $p>3$.
\end{thm}

In the subsequent work \cite{BeltranDuncanHickman}, Beltran, Duncan and Hickman
proved essentially optimal $L^p$ improving estimates for a local variant, with implications for essentially sharp sparse bounds for ${M^{\mathrm{hel}}}$.
It should also be interesting to explore versions of this work with restricted (fractal) dilation sets.

\subsection{The \texorpdfstring{Nevo--Thangavelu}{Nevo-Thangavelu} spherical means} \label{sec:NTh}
These integral operators were introduced by Nevo and Thangavelu in \cite{NevoThangavelu1997}. They
are considered interesting test cases for generalized Radon transforms associated with manifolds of higher co-dimension.
The concrete setup is as follows: Let $G$ be a two step nilpotent Lie group identified via the Baker-Campbell-Hausdorff formula with the Lie algebra ${{\mathbb{R}}}^{d}\times {{\mathbb{R}}}^m$; here the group law is given by
\begin{equation}\label{eq:grouplaw} (\underline x, \overline x)\cdot (\underline y,\overline y)=(\underline x+\underline y, \overline x+\overline y+\underline x^\intercal \vec J \underline y) ,\end{equation} with
$\underline x^\intercal \vec J \underline y=
(\underline x^\intercal J_1 \underline y,\dots \underline x^\intercal J_m \underline y)$ and the $J_i$ are $d\times d$ skew symmetric matrices. The case of $(d,m)$--M\'etivier groups is of particular interest; here it is assumed that $\sum_{i=1}^m \theta_i
J_i$ is invertible whenever $\theta\neq 0$. The Heisenberg group ${{\mathbb{H}}}_n$ arises for $d=2n$, $m=1$. For every $m>1$, M\'etivier groups exist but for given $m$ there are restrictive conditions on $d$ which are closely related to the Radon--Hurwitz numbers; in particular we must have $\log d {\gtrsim} m$. The explicit conditions can be found in \cite{AdamsLaxPhillips},
\cite{Kaplan1980}, see also the relevant work on generalized Radon transforms by Gressman \cite{gressman-RH}.

The parabolic dilations $\delta_t(x)=(t\underline x, t^2\overline x)$
are automorphisms of the group. We let $\mu$ be the normalized surface measure on the unit sphere in ${{\mathbb{R}}}^d$ and define the dilate $\mu_t$ by $\inn f{\mu_t}=\int f(t\underline y,0)\,d\mu$.
Then the Nevo--Thangavelu means are defined by the noncommutative convolution \[{{\mathscr{A}}}_t f(x)=f*\mu_t(x)= \int_{S^{d-1}} f(\underline x -t{\omega},\overline x -t\underline x^\intercal \vec J{\omega}) \,d\mu({\omega}).\]
Let ${{\mathfrak{M}}} f=\sup_{t>0} |{{\mathscr{A}}}_t f|$ be the associated maximal
function. Note that the support of the spherical measure is a co-dimension $m+1$ sphere in $G$. In the commutative case we have $\vec J=0$ and all results for ${{\mathfrak{M}}}$ can be obtained from the Euclidean case by a slicing argument. Thus
we assume $\vec J\neq 0$.

Following partial results in \cite{NevoThangavelu1997} the optimal range $p>\frac{d}{d-1}$, for $L^p$-boundedness of ${{\mathfrak{M}}}$ on M\'etivier groups was proved by M\"uller and the second author \cite{MuellerSeeger2004} for $d\ge 3$. Moreover, on the Heisenberg groups ${{\mathbb{H}}}_n$, $n\ge 2$, the result was independently obtained by
Narayanan and Thangavelu \cite{NarayananThangavelu2004}. The analysis in \cite{MuellerSeeger2004} relies on estimates for Fourier integral operators with folding canonical relations, applied to the spherical means. This is combined with a favorable behavior of the Schwartz kernel of $\frac{d}{dt}{{\mathscr{A}}}_t$ near the fold surface, due to the invariance of the horizontal subspace under the parabolic dilations.
The optimality is proved using a variant of the counterexample by Stein mentioned above. In contrast, the proof in \cite{NarayananThangavelu2004} is based on harmonic analysis on the Heisenberg groups. A third approach to the result on M\'etivier groups can be found in \cite{RoosSeegerSrivastava-IMRN2022}, using a rotational curvature property for $(x,t)\mapsto {{\mathscr{A}}}_t f(x)$ which results in a local smoothing phenomenon on $L^2$.
The idea goes back to a theorem by Joonil Kim \cite{KimJoonil-Annulus}, and is related to the $L^2$-analysis of Nikodym maximal functions in \cite{MockenhauptSeegerSogge1993}. We note that the square-function bound in \cite{NarayananThangavelu2004} also captures this local smoothing phenomenon.

It was asked in \cite{MuellerSeeger2004} whether the $L^p$ results for ${{\mathfrak{M}}}$ on M\'etivier groups extend to the general step two case. In this generality one can no longer use established results on boundedness of Fourier integral operators.
A partial $L^p$ result in this direction was obtained by
Liu--Yan \cite{LiuYan}, namely $L^p$ boundedness for
$d\ge 4$, $m=1$, and $p>\frac{d-1}{d-2}$. An optimal result was more recently proved by Jaehyeon Ryu and the second author \cite{RyuSeeger}:
\begin{thm}\label{thm:ryu-s} Let $G$ be any step two group with $d\ge 3$. Then ${{\mathfrak{M}}}$ is bounded on $L^p(G)$ if and only if $p>\frac{d}{d-1}$.
\end{thm}

In the proof of boundedness a linear algebra reduction is made to the situation where $J_1,\dots, J_m$ are linearly independent in the space of $d\times d$ skew symmetric matrices. After a dyadic frequency decomposition one interpolates between a weak type inequality for $L^1$ functions and an $L^2$ estimate; in both one takes advantage of the invariance of the horizontal space under the parabolic dilations.

Let us focus on the $L^2$ bound. Using essentially Plancherel's theorem one can reduce to the analysis
of a two-parameter family of oscillatory integral operators $T_{{\lambda},b}$ on $L^2({{\mathbb{R}}}^d)$. To define it let, for $x\in {{\mathbb{R}}}^d$, $x'=(x_2,\dots, x_d)$, and let $x'\mapsto g(x')$ be a real valued function satisfying
$g(0) = 1$, $g'(0) = 0$, and $g''(0)$ positive definite (reflecting the curvature of the sphere).
Let $P:x\mapsto x'$ be the projection to ${{\mathbb{R}}}^{d-1} $, let $S$
be a skew-symmetric $d\times d$ matrix of nonzero
rank, let $\rho(x',y_1) = y_1 + (x')^\intercal P Se_1$, let $a(x,t,y)$ be $C^\infty$ with small support where $t\approx 1$ and define the phase function $\psi\equiv\psi^S$ by
\begin{subequations}
\begin{equation}\psi(x,t,y)=y_1 \big(x_1-t g(\tfrac{x'-y'}{t})\big) +x^\intercal S\big( P^\intercal y'- t g(\tfrac{x'-y'}{t}) e_1\big) .\end{equation}
For ${\lambda}$ large and ${\lambda}^{-1}\le b\le 1$ define
\begin{equation} \label{eq:Tlab} T_{{\lambda},b} f(x,t)=\int e^{i{\lambda} \psi(x,t,y) }a(x,y,t) \chi_1 (\tfrac{ |t\rho(x',y_1)|}{b})f(y) \,dy.\end{equation}
\end{subequations} The quantity $\rho(x',y_1)$ is comparable in size with $\det \psi''_{xy}$ and should be thought of as an analogue of ``rotational curvature'' for the original problem; as $t\approx 1$ we localize in \eqref{eq:Tlab} to the set where $|\rho|\approx b$.
The crucial uniform space-time $L^2$ bound is
\begin{equation}\label{eq:Tlabbd} \|T_{{\lambda},b} \|_{L^2({{\mathbb{R}}}^d)\to L^2({{\mathbb{R}}}^{d+1})} {\lesssim} C_{\varepsilon} {\lambda}^{-\frac{d-1}{2}} b^{-\frac{d-2+{\varepsilon}}{2} } \end{equation} where the implicit constant depends on $S$ only in terms of the lower and upper bounds of the matrix norm $\|S\|$. We apply \eqref{eq:Tlab} for the choice $S=\sum_{i=1}^m \theta_i J_i$ where $|\theta|=1$.
One can view \eqref{eq:Tlabbd} as a phenomenon of “damping oscillatory integral operators” (see \cite{SoggeStein, SoggeStein1990}, \cite{PhongStein91} and many other such works in the literature).
Note that the amplitude $\chi_1 (\tfrac{ |t\rho(x',y_1)|}{b} )a$ loses regularity for small $b$; and it is crucial that $b$ is allowed in the full range where $b\ge {\lambda}^{-1}$. In \cite{RyuSeeger} inequality \eqref{eq:Tlabbd} is proved via several almost orthogonality arguments.

\subsubsection*{\texorpdfstring{\it The open case $d=2$.}{The open case}} The $L^2$ bound \eqref{eq:Tlabbd} is not enough for the estimation of the maximal operator when $d=2$, primarily because of the failure of a Sobolev embedding endpoint inequality for $L^\infty$. The following problem has attracted considerable attention:

\begin{problem}\label{opend=2}
Let $G$ be a step two group with $d=2$. Is ${{\mathfrak{M}}}$ bounded on $L^p(G)$ for $p>2$?
\end{problem}
This is currently open even in the model case of the Heisenberg group ${{\mathbb{H}}}_1$. We remark that for fixed $t$ the Nevo--Thangavelu means ${{\mathscr{A}}}_t$ on ${{\mathbb{H}}}_1$ (when acting on functions supported near the origin of the Heisenberg group) and the helical averages ${A^{\mathrm{hel}}}_t$ in ${{\mathbb{R}}}^3$ can be treated analogously; i.e. they both are Fourier integral operators of order $-1/2$ with a folding canonical relation and a curvature condition in the fibers above the fold surface (mentioned already in \S\ref{sec:helical} above). This is discussed in \cite{GreenleafSeeger1994} (see also the survey \cite{GreenleafSeegersurvey})
and unifies $L^p$-improving results in \cite{Oberlin-helix} and \cite{Secco}; more recently corresponding unifying sharp $L^p$-Sobolev estimates based on decoupling were proved in
\cite{PramanikSeeger2021}.

The analogy between the two families ${A^{\mathrm{hel}}}$ and ${{\mathcal{A}}}$ breaks down if we incorporate the $t$-variable. We view both as families of generalized Radon transforms acting on functions supported in a small neighborhood of $0$. In view of the above-mentioned curvature condition for the fold surface we now disregard
it and consider the geometry of the canonical relations associated with these Fourier integral operators away from the fold surface, i.e. microlocalized to parts which are graphs of conic symplectic diffeomorphisms.
There one uses the notion of {\it cinematic curvature}, introduced by Sogge in \cite{Sogge91} and further elaborated upon in \cite{ MockenhauptSeegerSogge1993}. The family
$\{{A^{\mathrm{hel}}}_t\}_{t\approx 1}$ of helical averages satisfies a {\it rank one} cinematic curvature condition which allows for an efficient use of decoupling techniques, leading in \cite{PramanikSeeger} to satisfactory estimates at least in the range $p>4$ ({\emph{cf.\ }} Theorem \ref{thm:helicalmax}). In contrast, the family
$\{{{\mathscr{A}}}_t\}_{t\approx 1}$ of Nevo--Thangavelu averages on ${{\mathbb{H}}}_1$ satisfies a {\it rank two} cinematic curvature condition which
does not provide a favorable numerology for an application of decoupling.

The subspace $L^p_{\mathrm{Hrad}} $ of Heisenberg-radial functions on ${{\mathbb{H}}}_1$ (i.e. $L^p$ functions of the form $f(x,u)=f_\circ (|\underline x|,u)$) is invariant under the action of ${{\mathscr{A}}}_t$. In order to gather evidence for an affirmative answer to the question in Problem \ref{opend=2}, Beltran, Guo, Hickman and the second author
\cite{BeltranGuoHickmanSeeger-Pisa2022}
showed that ${{\mathfrak{M}}}:L^p_{\mathrm{Hrad}}({{\mathbb{H}}}_1)\to L^p({{\mathbb{H}}}_1)$ is bounded for all $p>2$.
This can be reduced to a variable coefficient problem on maximal functions for certain weakly singular Radon transforms in ${{\mathbb{R}}}^2$. Despite the additional singularities one can use either techniques in \cite{MockenhauptSeegerSogge1992} or decoupling techniques in \cite{BeltranHickmanSogge} to obtain this result. Some simplifications and more on $L^p$-improving estimates on the space of Heisenberg-radial functions can be found in a subsequent paper by J. Lee and S. Lee \cite{LeeLee-Heisenberg-MA2023}.

\subsubsection*{\texorpdfstring{\it $L^p$-improving estimates}{The Lp improving estimates}}
We now return to the case $d\ge 3$ and discuss $L^p\to L^q$ bounds for the local maximal function ${{\mathfrak{M}}}_{\mathrm{loc}} f(x)= \sup_{1\le t\le 2} |{{\mathscr{A}}}_t f| $.
This problem was first investigated by Bagchi, Hait, Roncal and Thangavelu \cite{BagchiHaitRoncalThangavelu2021} for the Heisenberg groups ${{\mathbb{H}}}_n$ in order to obtain sparse domination results.
In this context the sparse forms are defined via dyadic cubes for the Heisenberg geometry; see also the more recent paper \cite{Conde-AlonsoDiPlinioParissisVempati} for sparse domination results in a very general geometric setting. Essentially optimal $L^p$ improving results (at least on the Heisenberg groups) were obtained in our work with Rajula Srivastava \cite{RoosSeegerSrivastava-IMRN2022} (see also \cite{RoosSeegerSrivastava-Studia2023} for variants with restricted dilation sets).

\begin{thm}[\cite{RoosSeegerSrivastava-IMRN2022}]\label{thm:mainLpLq} Let $G$ be a $(d,m)$--M\'etivier group, with $d>2$.
Let $(p^{-1}, q^{-1})$ belong to the interior of the quadrilateral with corners
\[ (0,0), \,(\tfrac{d-1}{d}, \tfrac{d-1}{d}), (\tfrac{d+m-1}{d+2m}, \tfrac{m+1}{d+2m}), (\tfrac{(d+m)(d+m-1)}{(d+m)^2+(d+m+1)m+1}, \tfrac{(m+1)(d+m-1)}{(d+m)^2+(d+m+1)m+1}).
\]
Then
${{\mathfrak{M}}}_{\mathrm{loc}}$ is bounded from $L^p(G)$ to $L^q(G)$.
\end{thm}
Various endpoint results are also obtained in \cite{RoosSeegerSrivastava-IMRN2022} (note that our notation differs from that one in \cite{RoosSeegerSrivastava-IMRN2022}). The basic strategy in the proof of theorem \ref{thm:mainLpLq} is as in the Euclidean problem. In the proof, it is crucial to observe a favorable cinematic curvature condition which holds independently of the codimension of the singular supports (in particular independent of $m$). We suspect that the result of Theorem \ref{thm:mainLpLq} is essentially sharp, but this has so far only been proved in \cite{RoosSeegerSrivastava-IMRN2022} for $m=1$, by a nonconventional Knapp type example.

\begin{problem} Let $G$ be a $(d,m)$--M\'etivier group. Is the result of Theorem \ref{thm:mainLpLq} sharp (up to endpoints) when $m\ge 2$?
\end{problem}

Next, dropping the M\'etivier condition, we can ask for $L^p$ improving properties for the local Nevo-Thangavelu maximal function on a general step two group. A slicing argument in \cite{RyuSeeger} shows that these cannot hold if $J_1,\dots, J_m$ are linearly dependent in the space of skew-symmetric matrices. Otherwise, the question is wide open.

\begin{problem} On a general $2$-step group determine essentially sharp $L^p\to L^q$ estimates for the maximal operator ${{\mathfrak{M}}}_{{\text{\rm loc}}}$.
\end{problem}

{{\it Remarks.}}
(i) In place of the Nevo--Thangavelu maximal function one can also consider the more traditional means for the codimension one {\it Kor\'anyi-spheres} $\{x:|\underline x|^4+ |\overline x|^2=t^4\}$ on the Heisenberg group, see \cite{Cowling1983}, \cite{Schmidt-Diplomarbeit}, \cite{GangulyThangavelu}. Rajula Srivastava \cite{SrivastavaR2024} proved essentially optimal $L^p\to L^q$-improving bounds for an associated local maximal function.

(ii) Joonil Kim \cite{KimJoonil2025} obtained $L^p$ estimates for another interesting spherical maximal function on ${{\mathbb{R}}}^{2n}\times {{\mathbb{R}}}$, working with a group law where $J_1$ is replaced by a not necessarily skew-symmetric matrix.

(iii) It should be very interesting to pursue the $L^p$ boundedness questions for variants of the Nevo--Thangavelu maximal operator in higher step groups.

\section{A fractal local smoothing conjecture}\label{sec:lsm}
Here we consider the half wave operators $e^{it\sqrt{-\Delta}}$ which are related to the spherical means by composing them with (pseudo-)differential operators of order $\frac{d-1}{2}$.
We discuss some $L^p$-Sobolev space-time estimates for solutions of the wave equation which are valid for larger $L^p$-Sobolev classes than the fixed time estimates; in view of the time integrals these have been labeled {\it local smoothing estimates.} Such inequalities had been early applied to simplify Bourgain's proof \cite{BourgainJdA1986} of the boundedness of the circular maximal operator in two dimensions (\cite{MockenhauptSeegerSogge1992}). Here we are interested in essentially sharp forms of this local smoothing phenomenon.
Recall Sogge's conjecture \cite{Sogge91} according to which one should have for $p>2$ the inequality
\begin{equation}\label{sogge-conjecture}
\Big(\int_1^2 \big
\| e^{i t \sqrt{-\Delta}} f \big \|_{p}^p\, d t\Big)^{\frac 1p} \leq C_{\varepsilon} \| f \|_{L^p_{ \sigma_p + \varepsilon}}, \, \,\, \sigma_p =\max\{0, \tfrac{d-1}2-\tfrac dp\}
\end{equation}
for all $\varepsilon >0$ and $2 < p <\infty$. This states that for $p>\frac{2d}{d-1}$ one gains a derivative of order $\frac 1p-{\varepsilon}$ when compared to the $L^p_s \to L^p$ fixed time result by Peral \cite{Peral1980} and Miyachi \cite{MiyachiWave}, with $s=(d-1)(\tfrac 12-\tfrac 1p)$.
To test the conjecture, its rather straightforward version on radial functions was verified early in \cite{MuellerSeeger1995}, see also \cite{ColzaniCominardiStempak} for an endpoint result.
The first positive result in the general case was proved in a groundbreaking paper by Wolff \cite{Wolff2000} who introduced the idea of decoupling inequalities for this problem (see also the subsequent papers \cite{LabaWolff} and \cite{GarrigosSchlagSeeger}). Sogge's conjecture was refined by Heo, Nazarov and the second author \cite{HeoNazarovSeeger-Acta} who conjectured \eqref{sogge-conjecture} with ${\varepsilon}=0$ for $p>\frac{2d}{d-1} $. This endpoint estimate was proved in \cite{HeoNazarovSeeger-Acta} for
$p>\frac{2(d-1)}{d-3}$, and $d\ge 4$ (no such endpoint bound is currently available in dimensions two and three).
Bourgain and Demeter \cite{BourgainDemeter2015} obtained \eqref{sogge-conjecture} for $p\ge \frac{2(d+1)}{d-1}$, via their optimal decoupling theorem for cones. Guth, Wang and Zhang \cite{GWZ} proved the full Sogge conjecture for $d=2$. In higher dimensions the currently best result is due to Gan, He, Li and Wu \cite{GanHeLiWu} who established \eqref{sogge-conjecture} in ranges that are not covered by the decoupling theorem.
We note that
the analogue of Sogge's conjecture in the $L^p_{\mathrm{rad}} (L^2_{\mathrm{sph}})$-category (which is implied by the $L^p$ version) was proved by M\"uller and the second author \cite{MuellerSeeger2001}.

In order to formulate a fractal version of Sogge's conjecture we consider an operator $P_j=\beta(2^{-j} |D|)$ with $\beta$ smooth and compactly supported away from $0$; thus $P_j$ localizes the Fourier transform of a function $f$ to an annulus $\{|\xi|\approx 2^j\}$. Let $\{t_\nu\}$ be any $2^{-j}$-separated subset of $[1,2]$. Then it follows from a Plancherel-P\'olya argument that
\eqref{sogge-conjecture} for $2\le p<\infty$ is equivalent to
\[\Big(\sum_\nu \big \|e^{it_\nu \sqrt{-\Delta}} P_j f \big\|_p^p\Big)^{1/p} {\lesssim} \|f\|_{L^p_s} \quad\text{
for $s>\max\big\{\tfrac 1p,(d-1)(\tfrac 12-\tfrac 1p) \big\}$}.
\]
For a fractal version, we are given a set $E\subset [1,2]$ and simply replace
$\{t_\nu\}$ with a $2^{-j} $-discretization $E(j)$, i.e. any set of $2^{-j}$-separated points in $E$. We are then interested in how the estimate in Sogge's conjecture improves. Since we will never be able to beat the threshold $(d-1)(\frac 12-\frac 1p)$ for the fixed time estimate such an improvement for $p>2$ can only occur when $(d-1)(\frac 12-\frac 1p)<\tfrac 1p$ i.e. when $p<\frac{2d}{d-1}$.
The behavior conjectured in \cite{BeltranRoosRutarSeeger} again features the Legendre--Assouad function $\nu_E^\sharp$
from \eqref{eq:LA}, \eqref{eq:ass-leg}.

\medskip
\begin{problem}\label{fractlsm}
For $p\ge 2$ let ${\kappa}_E(p)= \tfrac1p\nu_E^\sharp\big( (d-1) (\tfrac p2-1)\big)\,.$
Does the inequality
\begin{equation}\label{fractal-lsm-conj} \Big(\sum_{t\in E(j)}\big \|e^{it \sqrt{-\Delta} }P_j f\big \|_p^p\Big)^{1/p} \le C_{s,p} \|f\|_{L^p_s}, \quad s> {\kappa}_E(p), \end{equation}
hold for all $j>0$ and all $2^{-j}$-discretizations $E(j)$ of $E$?
\end{problem}

This would be sharp since in \cite{BeltranRoosRutarSeeger} it is shown that \eqref{fractal-lsm-conj} fails for some $f\in L^p$ if
$s<{\kappa}_E(p)$. Moreover, it is proved there that the conjecture holds for the class of {\it radial} $L^p_s$ functions.

\medskip

\begin{examples}
(i) If $E\subset[1,2]$ is quasi-Assouad regular with $\mathrm{dim}_{\mathrm{M}}\,E=\beta$ and $\mathrm{dim}_{\mathrm{qA}}\,E=\gamma$, then in
Problem \ref{fractlsm},
\[ {\kappa}_E(p) = \begin{cases} \big(1-\tfrac{\beta}{\gamma}\big) (d-1)(\tfrac 12-\tfrac 1p) + \tfrac{\beta}{p} \quad &\text{ if } \, 2 \leq p \leq 2+\frac{2{\gamma}}{d-1},
\\ (d-1)(\tfrac 12-\tfrac 1p)
\quad &\text{ if } \, p\ge 2+\frac{2{\gamma}}{d-1} \,.
\end{cases}
\]
In particular, for $E=[1, 2]$ one has ${\kappa}_E(p)= \max\{\frac 1p, (d-1)(\frac 12-\frac 1p) \} $ and thus one recovers Sogge's conjecture.

(ii) For $p\ge 2+\frac{2\gamma}{d-1}$ we always have
\[{\kappa}_E(p)=(d-1)(\tfrac 12-\tfrac 1p),\]
where $\mathrm{dim}_{\mathrm{qA}}\,E=\gamma$.
For general $E$, the behavior of ${\kappa}_E(p)$ for $2\le p\le 2+\frac{2{\gamma}}{d-1}$ can be fairly arbitrary (and can be characterized precisely using Proposition \ref{prop:whichhassleg}).

\end{examples}

\smallskip

The Legendre--Assouad function is used in \cite{BeltranRoosRutarSeeger} to formulate a more general off-diagonal $L^p\to L^q(\ell^q_{E(j)})$ conjecture.

\begin{conjecture} \label{general-lsm-conjecture}
Let $E \subset [1,2]$ and $1 < p \leq q < \infty$, $q > p'$. Then for every \[s> s_E(p,q):=\tfrac{d+1}{2}(\tfrac{1}{p}-\tfrac{1}{q}) + \tfrac{1}{q}\nu_E^\sharp \big( \tfrac{q(d-1)}{2}(1-\tfrac{1}{p} - \tfrac{1}{q}) \big)\] there exists a constant $C_{s,p,q} >0$ such that
\[\Big(\sum_{t\in E_j} \big\|e^{it\sqrt{-\Delta}} P_j
f\big \|^q_q\Big)^{\frac1q} \leq C_{s,p,q} 2^{js}
\|f\|_p
\]
holds for all $j\ge 1$ and all $2^{-j}$-discretizations $E_j$ of $E$.
\end{conjecture}
A proof of the conjecture would in particular imply a solution to Problem \ref{LAconj}. The case $p=2$ (verified in \cite{BeltranRoosRutarSeeger}) amounts to a fractal version of a Strichartz inequality. It would be very interesting to obtain resolutions of this conjecture for other pairs of exponents.

\section*{Acknowledgments} We would like to thank David Beltran and Po Lam Yung for helpful and encouraging comments.
A.S. is grateful for the opportunity to present parts of this survey in an invited analysis section talk at the International Congress of Mathematicians 2026 in Philadelphia, PA.
J.R. was supported in part by a grant from the Simons Foundation and by National Science Foundation grant DMS-2154835. A.S. was supported in part by National Science Foundation grants DMS-2348797, DMS-2529637.

\end{document}